\documentclass[dvips]{arxbj}
\usepackage{cursive}   %{uplgreek}

\aid{0}
\volume{16}
\issue{2}
\pubyear{2010}
\firstpage{301}
\lastpage{330}
\doi{10.3150/09-BEJ201}

\makeatletter

\newtheorem{lemma}{Lemma}
\newremark{remark}{Remark}
\newtheorem{corollary}[lemma]{Corollary}
\newtheorem{theorem}[lemma]{Theorem}
\makeatother

\begin{document}
\begin{frontmatter}

\title{A new method for obtaining sharp compound Poisson
approximation error estimates for sums of locally dependent random variables}
\runtitle{Sharp compound Poisson
approximation error estimates}

\begin{aug}
\author[1]{\fnms{Michael V.} \snm{Boutsikas}\thanksref{1}\ead[label=e1]{mbouts@unipi.gr}\corref{}}
\and
\author[2]{\fnms{Eutichia} \snm{Vaggelatou}\thanksref{2}\ead[label=e2]{evagel@math.uoa.gr}}
\runauthor{M.V. Boutsikas and E. Vaggelatou}
\address[1]{Department of Statistics and Insurance
Science, University of Piraeus,
 Karaoli and Dimitriou Str. 80, Piraeus 18534, Greece. \printead{e1}}
\address[2]{Section of Statistics and Operations Research,
 Department of Mathematics, University of Athens,
Panepistemiopolis, Athens 15784, Greece. \printead{e2}}
\end{aug}

% HISTORY:
\received{\smonth{3} \syear{2008}}
\revised{\smonth{11} \syear{2008}}

% ABSTRACT
\begin{abstract}
Let $X_{1},X_{2},\ldots,X_{n}$ be a sequence of independent or locally
dependent random variables taking values in $\mathbb{Z}_{+}$. In this paper,
we derive sharp bounds, via a new probabilistic method, for the total
variation distance between the distribution of the sum $\sum_{i=1}^{n}X_{i}$
and an appropriate Poisson or compound Poisson distribution. These bounds
include a factor which depends on the smoothness of the approximating
Poisson or compound Poisson distribution. This ``smoothness
factor'' is of order $\mathrm{O}(\sigma ^{-2})$, according to a
heuristic argument, where~$\sigma ^{2}$ denotes the variance of the
approximating distribution. In this way, we offer sharp error estimates for
a large range of values of the parameters. Finally, specific examples
concerning appearances of rare runs in sequences of Bernoulli trials are
presented by way of illustration.
\end{abstract}

% KEYWORDS

\end{frontmatter}

%s1 ###
\section{Introduction and overview}\label{sec1}

Let $X_{1},X_{2},\ldots ,X_{n}$ be a sequence of independent or locally
dependent random variables which take values in $\mathbb{Z}_{+}$. If $%
X_{1},X_{2},\ldots,X_{n}$ rarely differ from zero (that is, $P(X_{i}\neq 0)\approx
0 $), then it is well known that the distribution of their sum can be
efficiently approximated by an appropriate Poisson or compound Poisson
distribution. This situation appears in a great number of applications
involving locally dependent and rare events, such as risk theory, extreme
value theory, reliability theory, run and scan statistics, graph theory and
biomolecular sequence analysis.

The main method used so far for establishing effective Poisson or compound
Poisson approximation results in the case of independent or dependent random
variables is the much acclaimed Stein--Chen method (see, for example, Barbour, Holst
and Janson (\citeyear{1992BarbourJan}), Barbour and Chryssaphinou (\citeyear{2001Barbour}),
Barbour and Chen (\citeyear{2005Barbour})
and the references therein). Another method for independent random variables
is Kerstan's method (see Roos (\citeyear{2003Roos}) and the references therein).

In the recent years, an alternative methodology has been developed in a
series of papers concerning compound Poisson approximation for sums or
processes of dependent random variables, employing probabilistic techniques,
that is, properties of certain probability metrics, stochastic orders and
coupling techniques  (see
Boutsikas and Koutras (\citeyear{2000Boutsikas}, \citeyear{2001Boutsikas}),
Boutsikas and Vaggelatou (\citeyear{2002Boutsikas}),
Boutsikas (\citeyear{2006Boutsikas})). In this series of papers, the error estimates are, under
analogous assumptions, of almost the same nature and the same order as the
error estimates developed by the Stein--Chen method. The main shortcoming
of these bounds, though, is that they do not incorporate any so-called
``magic factor'' (however, in the process
approximation case treated in Boutsikas (2006), such a factor cannot be
present). This factor, also known as a \textit{Stein factor}, appears in
approximation error estimates obtained through the Stein--Chen method and
decreases as the parameter of the Poisson distribution increases.

The purpose of this work is to derive sharp error bounds for the total
variation distance between the distribution of the sum of integer-valued
random variables and an appropriate Poisson or compound Poisson
distribution. Specifically, by assuming that the random variables $%
X_{1},X_{2},\ldots,X_{n}$ are locally dependent (in the strict sense of $k$%
-dependence), we derive bounds similar in nature to those obtained by
the Stein--Chen method that include a factor analogous to a Stein factor. This
factor is better/smaller than the associated Stein factors, thereby
offering (for a large range of the values of the parameters) sharper bounds
than relative ones derived via the Stein--Chen method. This factor is just
the $L_{1}$-norm, $\Vert \Delta ^{2}f\Vert _{1}$, of the second
difference of the probability distribution function $f$ of the approximating
Poisson or compound Poisson distribution. It decreases as $f$ becomes
smoother, which, in our case, usually happens when the variance of the
distribution corresponding to $f$ increases. Hence, we shall often refer to
this factor as the \textit{smoothness factor}. The methodology we employ is
based on a modification of Lindeberg's method, along with the coupling
inequality of Lemma \ref{lemma3} and the smoothing inequality (which produces the
aforementioned smoothness factor)  of Lemma \ref{Btv}.

It is worth pointing out an undesired effect of our treatment, which is an
additional term in the proposed bounds that does not appear in Stein--Chen
bounds. This term becomes large for a certain range of values of the
parameters, but, as we explain in Remark \ref{r3} of Section \ref{sec3}, it can be
substantially reduced if we possess a simple and effective upper bound for $%
\Vert \Delta ^{2}f\Vert _{1}$. Nevertheless, this term is generally negligible, especially for small or moderate values of
$\lambda $, where $\lambda $ is the parameter of the approximating Poisson
distribution.

It is worth stressing that the error estimates presented in this work have
the same optimal order as other bounds obtained through the Stein--Chen
method. In fact, bounds derived using the latter method contain an
additional $\log \lambda $ term or, worse, an $\mathrm{e}^{\lambda }$ term for
certain ranges of the parameters (see Barbour, Chen and Loh (\citeyear{1992Barbour}), Barbour
and Utev (\citeyear{1999Barbour}), Barbour and Xia (\citeyear{2000Barbour}) or Barbour and Chryssaphinou (\citeyear{2001Barbour})
and the references therein). On the other hand, our bounds do not include
such terms and they incorporate a better and more natural factor which we
conjecture to be optimal.

The paper is organized as follows. In Section \ref{sec2}, we present some already
known, as well as new, auxiliary lemmas which concern probability metrics
and coupling techniques. These lemmas will be used for the derivation of our
main results. In Sections \ref{sec3} and \ref{sec4}, we present our main results, that is, bounds concerning Poisson and compound Poisson approximation for sums of
\textit{independent} and $k$-\textit{dependent} random variables,
respectively. Finally, in Section \ref{sec5}, in order to illustrate the
applicability and  effectiveness of our main results, we present a simple
example of an application which concerns the distribution of appearances of
rare runs in sequences of independent and identically distributed (i.i.d.)
trials.

%s2 ###
\section{Preliminary results}\label{sec2}

Throughout this paper, the abbreviations c.d.f.~and p.d.f.~will stand for
the cumulative distribution function and probability density function,
respectively. In addition, $\mathcal{L} X$ or $\mathcal{L}(X)$ will denote
the distribution of a random variable $X$ and the notation $X\sim $ $G$ will
imply that $X$ follows the distribution $G$. Moreover, we shall write $%
\mathit{Po}(\lambda )$ to denote the Poissson distribution with mean $\lambda $ and $%
\mathit{CP}(\lambda ,F)$ to denote the compound Poisson distribution with Poisson
parameter $\lambda $ and compounding distribution $F$. In other words, $%
\mathit{CP}(\lambda ,F)$ is the distribution of the random sum $\sum_{i=1}^{N}Z_{i}$,
where $N\sim \mathit{Po}(\lambda )$ and $Z_{i}$ are i.i.d.~random variables with
c.d.f.~$F$ which are also independent of $N$. For two functions $f$ and $g$,
the following standard notation will be used:
\begin{eqnarray*}
f(t)\sim g(t)  \mbox{ as } t\rightarrow t_{0} \mbox{ if
}\lim_{t\rightarrow t_{0}}\frac{f(t)}{g(t)}=1 ;  \qquad f(t)=\mathrm{O}(g(t))  %
\mbox{ if  }\frac{f(t)}{g(t)}\mbox{ is bounded.}
\end{eqnarray*}%
Moreover, whenever dependence or independence of some random variables is
mentioned, it will be immediately assumed that they are defined on the same
probability space. Finally, $\lfloor x\rfloor $ denotes the
integer part of $x$ and we will assume that $\sum_{i=a}^{b}x_{i}=0$ when $%
a>b $.

%s2.1 ###
\subsection{Probability metrics and smoothness factors}\label{sec2.1}

In order to quantify the quality of a distribution approximation, the
\textit{total variation distance} and \textit{Zolotarev's ideal metric} of
order $2$ will be used. Since the results of this paper concern discrete
distributions, it suffices to consider only the discrete versions of the
aforementioned probability metrics.

The \textit{total variation distance} between the distributions $\mathcal{L}%
X $ and $\mathcal{L}Y$ of two random variables $X$ and $Y$ is defined by%
\begin{eqnarray*}
d_{\mathrm{TV}}( \mathcal{L}X,\mathcal{L}Y) :=\sup_{A\subseteq \mathbb{Z}%
}\vert P( X\in A) -P( Y\in A) \vert =\frac{1%
}{2}\sum_{k=-\infty }^{\infty }\vert P(X=k)-P(Y=k)\vert ,
\end{eqnarray*}
whereas the \textit{total variation distance of order }$2$ or \textit{%
Zolotarev's ideal metric of order} $2$ (Zolotarev (\citeyear{1983Zolotarev})) is defined by
\begin{eqnarray*}
\mathbf{\zeta }_{2}(\mathcal{L}X,\mathcal{L}Y)=\int_{-\infty }^{\infty
}\vert E(X-t)_{+}-E(Y-t)_{+}\vert\, \mathrm{d}t=\sum_{k=-\infty }^{\infty
}\Biggl\vert \sum_{u=k}^{\infty }\bigl( F_{X}(u)-F_{Y}(u)\bigr) \Biggr\vert ,
\end{eqnarray*}%
where, as usual, $F_{X}$ denotes the c.d.f.~of the random variable $X$.
Throughout, whenever a $\mathbf{\zeta }_{2}(\mathcal{L}X,\mathcal{L}Y)$
distance appears, it will be implicitly assumed that $X,Y$ \textit{possess
finite first and second moments} and that $E(X)=E(Y)$. For a comprehensive
exposition on probability metrics and their properties, the interested
reader may consult Rachev (\citeyear{1991Rachev}) and the references therein.

Next, we denote by $\Delta ^{k}f$ the $k$th order (backward) difference
operator over a function $f\dvtx \mathbb{Z\rightarrow R}$, that is, $\Delta
f(i)=f(i)-f(i-1)$ and $\Delta ^{k}=\Delta (\Delta ^{k-1}f)$, $k=1,2,\ldots
(\Delta ^{0}f=f)$. The smoothness factor mentioned in the
\hyperref[sec1]{Introduction} emerges from the following lemma. Analogous results concerning random
variables with a Lebesgue density have been used in the past in order to
obtain Berry--Esseen-type results (see Senatov (\citeyear{1980Senatov}), Rachev (\citeyear{1991Rachev}) and the
references therein).

\begin{lemma}
\label{Btv}If $X,Y,Z$ are integer-valued random variables (with finite first
and second moments) such that $E(X)=E(Y)$ and $Z$ is independent of $X,Y$,
then%
\begin{eqnarray*}
d_{\mathrm{TV}}\bigl(\mathcal{L}(X+Z),\mathcal{L}(Y+Z)\bigr)\leq \tfrac{1}{2}\Vert \Delta
^{2}f_{Z}\Vert _{1} \mathbf{\zeta }_{2}(\mathcal{L}X,\mathcal{L}Y),
\end{eqnarray*}
where $f_{Z}$ is the p.d.f.~of $Z$ and $\Vert \Delta
^{2}f_{Z}\Vert _{1}:=\sum_{z\in \mathbb{Z}}|\Delta ^{2}f_{Z}(z)|$.
\end{lemma}

\begin{pf}
For any functions $a,b\dvtx \mathbb{Z}\rightarrow \mathbb{R}$ and $c,d\in \mathbb{%
Z}$, we have (second-order Abel summation formula)
%
%e1 ###
\begin{equation}\label{eq48}
\sum_{z=c}^{d}b_{z-2} \Delta ^{2}a_{z}=\sum_{z=c}^{d}a_{z} \Delta
^{2}b_{z}+b_{d-1} \Delta a_{d}-a_{d} \Delta b_{d}+a_{c-1} \Delta
b_{c-1}-b_{c-2} \Delta a_{c-1}.
\end{equation}
Denote by $f_{W}$ the p.d.f. of any discrete random variable $W$. If,
for fixed $k,$ we now choose
\begin{eqnarray*}
a_{z}=f_{Z}(z),\qquad   b_{z}=\sum_{i=-\infty }^{z+1}\bigl(
R_{X}(k-i)-R_{Y}(k-i)\bigr) ,
\end{eqnarray*}
where $R_{X}(k-z)=\sum_{i=-\infty }^{z-1}f_{X}(k-i)$, and then take $%
c\rightarrow -\infty $, $d\rightarrow \infty $, identity (\ref{eq48}) leads
to
%
%e2 ###
\begin{equation} \label{eq3}
\sum_{z=-\infty }^{\infty }\sum_{i=-\infty }^{z-1}\bigl(
R_{X}(k-i)-R_{Y}(k-i)\bigr) \Delta ^{2}f_{Z}(z)=\sum_{z=-\infty }^{\infty
}\bigl(f_{X}(k-z)-f_{Y}(k-z)\bigr)f_{Z}(z)
\end{equation}
since all quantities $b_{z},a_{z},\Delta a_{z},\Delta b_{z} $vanish as $%
z\rightarrow \infty $ or $z\rightarrow -\infty $. Using (\ref{eq3}), we get
\begin{eqnarray*}
d_{\mathrm{TV}}\bigl(\mathcal{L}(X+Z),\mathcal{L}(Y+Z)\bigr) &=&\frac{1}{2}\sum_{k=-\infty
}^{\infty }\Biggl\vert \sum_{z=-\infty }^{\infty }\bigl(
f_{X}(k-z)-f_{Y}(k-z)\bigr) f_{Z}(z)\Biggr\vert
\\
&=&\frac{1}{2}\sum_{k=-\infty }^{\infty }\Biggl\vert \sum_{z=-\infty
}^{\infty }\Delta ^{2}f_{Z}(z)\sum_{i=-\infty }^{z-1}\bigl(
R_{X}(k-i)-R_{Y}(k-i)\bigr) \Biggr\vert
\\
&\leq &\frac{1}{2}\sum_{z=-\infty }^{\infty }\vert \Delta
^{2}f_{Z}(z)\vert \sum_{k=-\infty }^{+\infty }\Biggl\vert
\sum_{i=-\infty }^{z-1}\bigl( R_{X}(k-i)-R_{Y}(k-i)\bigr) \Biggr\vert .
\end{eqnarray*}
Finally, setting $s:=k-z+1$ and $u:=k-i$ in the second and third summation
above yields
\begin{eqnarray*}
\hspace{32pt}d_{\mathrm{TV}}\bigl(\mathcal{L}(X+Z),\mathcal{L}(Y+Z)\bigr) &\leq &\frac{1}{2}\sum_{z=-\infty
}^{\infty }\vert \Delta ^{2}f_{Z}(z)\vert \sum_{s=-\infty
}^{\infty }\Biggl\vert \sum_{u=s}^{\infty }\bigl( R_{X}(u)-R_{Y}(u)\bigr)
\Biggr\vert
\\
&=&\frac{1}{2}\Vert \Delta ^{2}f_{Z}\Vert _{1} \mathbf{\zeta }%
_{2}(X,Y).\hspace{132pt}\qed
\end{eqnarray*}
\noqed\end{pf}

If $Z=0$ and $E(X)=E(Y)$, then a simple consequence of the above result is
the inequality
%e3 ###
\begin{equation}\label{eq18}
d_{\mathrm{TV}}(\mathcal{L}X,\mathcal{L}Y)\leq \tfrac{1}{2}\Vert \Delta
^{2}f_{0}\Vert _{1} \mathbf{\zeta }_{2}(\mathcal{L}X,\mathcal{L}Y)=2
\mathbf{\zeta }_{2}(\mathcal{L}X,\mathcal{L}Y),
\end{equation}
where $f_{0}:=f_{Z}$ when $Z=0$. If $Z$ follows a Poisson distribution with
parameter $\lambda $, then we can find the explicit value of $\Vert
\Delta ^{2}f_{Z}\Vert _{1}$ and its asymptotic behavior. In the
sequel, we shall write $f_{\mathit{Po}(\lambda )}$ instead of $f_{Z}$ when $Z$ $\sim
\mathit{Po}(\lambda )$. As we will see below, it is convenient to first find the $%
L_{\infty }$-norm, $\Vert \Delta f_{\mathit{Po}(\lambda )}\Vert _{\infty }$%
, and then to investigate its relation with the norm $\Vert \Delta
^{2}f_{\mathit{Po}(\lambda )}\Vert _{1}$.

\begin{lemma}
\label{prop1}If $f_{\mathit{Po}(\lambda )}$ denotes the probability distribution
function of the Poisson distribution with parameter $\lambda $, then
\begin{eqnarray*}
\bigl\Vert \Delta f_{\mathit{Po}(\lambda )}\bigr\Vert _{\infty }=\sup_{k\in \mathbb{Z%
}_{+}}\bigl\vert f_{\mathit{Po}(\lambda )}(k)-f_{\mathit{Po}(\lambda )}( k-1)
\bigr\vert =\mathrm{e}^{-\lambda }\frac{\lambda ^{k_{\lambda }}}{k_{\lambda }!}%
\biggl( 1-\frac{k_{\lambda }}{\lambda }\biggr) ,
\end{eqnarray*}
where $k_{\lambda }:=\lfloor \lambda -\sqrt{\lambda +1/4}%
+1/2\rfloor $ for all $\lambda >0$. In particular, $\Vert \Delta
f_{\mathit{Po}(\lambda )}\Vert _{\infty }=\mathrm{e}^{-\lambda }$ for $\lambda \leq 2$.
Furthermore,
\begin{eqnarray*}
\bigl\Vert \Delta f_{\mathit{Po}(\lambda )}\bigr\Vert _{\infty }\sim \frac{1}{%
\lambda \sqrt{2\curpi \mathrm{e}}} \qquad  \mbox{as }\lambda \rightarrow \infty .
\end{eqnarray*}
\end{lemma}

\begin{pf}
It can be easily verified that $\Delta f_{\mathit{Po}(\lambda )}(k)=\mathrm{e}^{-\lambda }%
\frac{\lambda ^{k}}{k!}( 1-\frac{k}{\lambda })$, $k\in
\{0,1,2,\ldots\}$, while $\Delta f_{\mathit{Po}(\lambda )}(k)=0$ for $k<0$, and also that
\begin{eqnarray*}
\Delta ^{2}f_{\mathit{Po}(\lambda )}(k)=\mathrm{e}^{-\lambda }\frac{\lambda ^{k}}{k!}\biggl( 1+%
\frac{k(k-1)}{\lambda ^{2}}-2\frac{k}{\lambda }\biggr) ,\qquad  k\in \{0,1,2,\ldots\},
\end{eqnarray*}
while $\Delta ^{2}f_{\mathit{Po}(\lambda )}(k)=0$ for $k<0$. Define $h\dvtx \mathbb{R}%
_{+}\rightarrow \mathbb{R} $ such that $h( x) =\Delta
^{2}f_{\mathit{Po}(\lambda )}(x)$ (that is, the extension of $\Delta ^{2}f_{\mathit{Po}(\lambda
)} $ over $\mathbb{R}_{+})$, where  $x!$ now denotes the Gamma function $%
\Gamma (1+x)$. It is easy to verify that $h$ is positive when $0\leq x\leq
\rho _{1} $, negative when $\rho _{1}\leq x\leq \rho _{2}$ and positive
again when $x\geq \rho _{2}$, where $\rho _{1}=\rho _{1}(\lambda )=\lambda -%
\sqrt{\lambda +1/4}+1/2$ and $\rho _{2}=\rho _{2}(\lambda )=\lambda +\sqrt{%
\lambda +1/4}+1/2$ are the two roots of the equation $h( x) =0$ ($%
0<\rho _{1}<\rho _{2}$). Since $h$ is an extension of $\Delta
^{2}f_{\mathit{Po}(\lambda )}$, we deduce that $\Delta ^{2}f_{\mathit{Po}(\lambda )}(k)\geq 0$
when $0\leq k\leq \rho _{1}$, $\Delta ^{2}f_{\mathit{Po}(\lambda )}(k)\leq 0$ when $%
\rho _{1}\leq k\leq \rho _{2}$ and $\Delta ^{2}f_{\mathit{Po}(\lambda )}(k)\geq 0$
when $k\geq \rho _{2}$. This implies that $0=\Delta f_{\mathit{Po}(\lambda )}(-1)\leq
\Delta f_{\mathit{Po}(\lambda )}(0)\leq \cdots\leq \Delta f_{\mathit{Po}(\lambda )}(\lfloor
\rho _{1}\rfloor )$, while $\Delta f_{\mathit{Po}(\lambda )}(\lfloor \rho
_{1}\rfloor )\geq \Delta f_{\mathit{Po}(\lambda )}(\lfloor \rho
_{1}\rfloor +1)\geq \cdots\geq \Delta f_{\mathit{Po}(\lambda )}(\lfloor \rho
_{2}\rfloor )$  and $\Delta f_{\mathit{Po}(\lambda )}(\lfloor \rho
_{2}\rfloor )\leq \Delta f_{\mathit{Po}(\lambda )}(\lfloor \rho
_{2}\rfloor +1)\leq \cdots$. Hence, $|\Delta f_{\mathit{Po}(\lambda )}(k)|$ must
be maximized at $\lfloor \rho _{1}\rfloor $ or $\lfloor \rho
_{2}\rfloor $ (since $\Delta f_{\mathit{Po}(\lambda )}(k)\rightarrow 0$ as $%
k\rightarrow \infty $). In order to verify that it is maximized at $%
k_{\lambda }=\lfloor \rho _{1}(\lambda )\rfloor $, we shall prove
that $g_{1}(\lambda )>g_{2}(\lambda ) $ for all $\lambda >0$ where $%
g_{1}(\lambda )=\lambda |\Delta f_{\mathit{Po}(\lambda )}(\lfloor \rho
_{1}(\lambda )\rfloor )|$ and $ g_{2}(\lambda )=\lambda |\Delta
f_{\mathit{Po}(\lambda )}(\lfloor \rho _{2}(\lambda )\rfloor )|$, that is,
\begin{eqnarray*}
g_{1}(\lambda )&=&\lambda \mathrm{e}^{-\lambda }\frac{\lambda ^{\lfloor \rho
_{1}(\lambda )\rfloor }}{\lfloor \rho _{1}(\lambda )\rfloor
!}\biggl( 1-\frac{\lfloor \rho _{1}(\lambda )\rfloor }{\lambda }%
\biggr) >g_{2}(\lambda )
\\
&=&-\lambda \mathrm{e}^{-\lambda }\frac{\lambda ^{\lfloor
\rho _{2}(\lambda )\rfloor }}{\lfloor \rho _{2}(\lambda
)\rfloor !}\biggl( 1-\frac{\lfloor \rho _{2}(\lambda
)\rfloor }{\lambda }\biggr) ,\qquad \lambda >0.
\end{eqnarray*}
For every $k\in \{0,1,\ldots\},\varepsilon \in \lbrack 0,1)$, we have $%
\lfloor \rho _{1}(k+\varepsilon +\sqrt{k+\varepsilon })\rfloor
=\lfloor k+\varepsilon \rfloor =k$. Therefore, $\lfloor \rho
_{1}(\lambda )\rfloor =k$ for every $\lambda \in \lbrack k+\sqrt{k}%
,k+1+\sqrt{k+1})$. Hence, in this interval, the function $g_{1}(\lambda )$
is equal to $\lambda \mathrm{e}^{-\lambda }\frac{\lambda ^{k}}{k!}( 1-\frac{k}{%
\lambda }) $,  differentiable (except at $k+\sqrt{k})$ and concave,
and $g_{1}^{\prime }(\lambda )=0$ at $\lambda =a(k)=k+1/2+\sqrt{k+1/4}$.
Moreover, $g_{1}(\lambda )\rightarrow g_{1}(k+1+\sqrt{k+1})$ as $\lambda
\rightarrow k+1+\sqrt{k+1}$ and thus $g_{1}(\lambda )$ is continuous for
every $\lambda >0$. Therefore, $g_{1}(\lambda )\geq g_{1}(k+\sqrt{k})$ for
every $\lambda \in \lbrack a(k-1),a(k)]$, $k\in \{1,2,\ldots\}$. Using the
upper bound of Stirling's approximation ($k!\leq k^{k}\mathrm{e}^{-k}\sqrt{2\curpi k}\mathrm{e}^{%
{1/(12k)}}$) and the elementary inequality $\log
(1+x)>x-x^{2}/2+x^{3}/3-x^{4}/4, x>0$, we get
\begin{eqnarray*}
g_{1}\bigl(k+\sqrt{k}\bigr)=\mathrm{e}^{-(k+\sqrt{k})}\frac{(k+\sqrt{k})^{k}}{k!}\sqrt{k}\geq
\frac{\mathrm{e}^{-\sqrt{k}-{1/(12k)}}}{\sqrt{2\curpi }}\mathrm{e}^{k\log (1+{1/\sqrt{k}%
})}>\frac{\mathrm{e}^{{1/(3\sqrt{k})}-{1/(3k)}}}{\sqrt{2\curpi \mathrm{e}}}\geq \frac{1}{%
\sqrt{2\curpi \mathrm{e}}}
\end{eqnarray*}
for every $k\geq 1$. Therefore, $g_{1}(\lambda )>\frac{1}{\sqrt{2\curpi \mathrm{e}}}$
for every $\lambda \in \bigcup _{k\geq 1}[a(k-1),a(k)]=[1,\infty )$.

 Similarly, for every $k\in \{1,2,\ldots\}, \varepsilon \in \lbrack 0,1)$,
we have $\lfloor \rho _{2}(k+\varepsilon -\sqrt{k+\varepsilon })\rfloor
=\lfloor k+\varepsilon \rfloor =k$. Therefore, $\lfloor \rho
_{2}(\lambda )\rfloor =k$ for every $\lambda \in \lbrack k-\sqrt{k}%
,k+1-\sqrt{k+1})$. Moreover, in this interval, the function $g_{2}(\lambda )$
is equal to $\lambda \mathrm{e}^{-\lambda }\frac{\lambda ^{k}}{k!}( \frac{k}{%
\lambda }-1) $,  differentiable (except at $k-\sqrt{k})$ and
concave, and $g_{2}^{\prime }(\lambda )=0$ at $\lambda =k+1/2-\sqrt{k+1/4}$ (%
$g_{2}(\lambda )$ is also continuous for every $\lambda >0$). Therefore, $%
g_{2}(\lambda )\leq g_{2}(k+1/2-\sqrt{k+1/4})$ for every $\lambda \in
\lbrack k-\sqrt{k},k+1-\sqrt{k+1}], k\in \{1,2,\ldots\}$. Using the lower bound
of Stirling's approximation $(k!\geq k^{k}\mathrm{e}^{-k}\sqrt{2\curpi k})$ and the
elementary inequality $\log (1+x)<x-x^{2}/2,x\in (-1,0),$ we get, for $k\geq 1$,
\begin{eqnarray*}
g_{2}\bigl(k+{1/2}-\sqrt{k+{1/4}}\bigr) &\leq &\frac{\sqrt{k+{1/4}}%
-{1/2}}{\sqrt{2\curpi k}}\mathrm{e}^{-{1/2}+\sqrt{k+{1/4}}+k\log (1+%
{({1/2}-\sqrt{k+{1/4}})/k})}
\\
&<&\frac{\sqrt{k+{1/4}}-{1/2}}{\sqrt{k}}\frac{\mathrm{e}^{{(\sqrt{k+%
{1/4}}-{1/2})/2k}}}{\sqrt{2\curpi \mathrm{e}}}<\frac{1}{\sqrt{2\curpi \mathrm{e}}}.
\end{eqnarray*}
Therefore, $g_{2}(\lambda )<\frac{1}{\sqrt{2\curpi \mathrm{e}}}$ for every $\lambda \geq
0$. Hence, $g_{2}(\lambda )<\frac{1}{\sqrt{2\curpi \mathrm{e}}}<g_{1}(\lambda )$ for
every $\lambda \geq 1$. It now remains to show that $g_{2}(\lambda
)<g_{1}(\lambda )$ for every $0<\lambda <1$. This is easily verified since $%
g_{1}(\lambda )=\lambda \mathrm{e}^{-\lambda }$ for $\lambda <2$, while $g_{2}(\lambda
)=\mathrm{e}^{-\lambda }\lambda ( 1-\lambda ) $ for $\lambda \in \lbrack
0,2-\sqrt{2})$ and $g_{2}(\lambda )=\mathrm{e}^{-\lambda }\lambda ^{2}( 1-\frac{%
\lambda }{2}) $ for $\lambda \in \lbrack 2-\sqrt{2},3-\sqrt{3})$.%

Finally, from $(1+y)^{k_{\lambda }}=\mathrm{e}^{k_{\lambda }\log
(1+y)}=\mathrm{e}^{k_{\lambda }(y-{y^{2}/2}+\mathrm{o}(y^{2}))}$ with $y=(\lambda
-k_{\lambda })/k_{\lambda}$, we get
\begin{eqnarray*}
\mathrm{e}^{k_{\lambda }-\lambda }( \lambda /k_{\lambda }) ^{k_{\lambda
}}\rightarrow \mathrm{e}^{-1/2}\qquad\mbox{as  }\lambda \rightarrow \infty .
\end{eqnarray*}
From this fact and Stirling's formula, we get, as $\lambda \rightarrow \infty
$, that
\begin{eqnarray*}
\hspace{40pt}\lambda \Delta f_{\mathit{Po}(\lambda )}(k_{\lambda })=\lambda \mathrm{e}^{-\lambda }\frac{%
\lambda ^{k_{\lambda }}}{k_{\lambda }!}\biggl( 1-\frac{k_{\lambda }}{\lambda }%
\biggr) \sim \mathrm{e}^{k_{\lambda }-\lambda }\biggl( \frac{\lambda }{k_{\lambda }}%
\biggr) ^{k_{\lambda }}\frac{\lambda -k_{\lambda }}{\sqrt{2\curpi k_{\lambda }}}%
\rightarrow \frac{1}{\sqrt{2\curpi \mathrm{e}}}.\hspace{27pt}\qed
\end{eqnarray*}
\noqed\end{pf}

In the next lemma, we find the explicit value of $\Vert \Delta
^{2}f_{\mathit{Po}(\lambda )}\Vert _{1}$ and a convenient upper bound in terms
of $\Vert \Delta f_{\mathit{Po}(\lambda )}\Vert _{\infty }$.

\begin{lemma}
\label{np}If $f_{\mathit{Po}(\lambda )}$ denotes the p.d.f. of the Poisson
distribution with parameter $\lambda $, then
\begin{eqnarray*}
\bigl\Vert \Delta ^{2}f_{\mathit{Po}(\lambda )}\bigr\Vert _{1}=\sum_{z=0}^{+\infty
}\bigl|\Delta ^{2}f_{\mathit{Po}(\lambda )}(z)\bigr|=2\mathrm{e}^{-\lambda }\biggl( \frac{\lambda
^{k_{\lambda }-1}( \lambda -k_{\lambda }) }{k_{\lambda }!}-\frac{%
\lambda ^{u_{\lambda }-1}( \lambda -u_{\lambda }) }{u_{\lambda }!}%
\biggr) ,
\end{eqnarray*}
where $k_{\lambda }:=\lfloor \lambda -\sqrt{\lambda +1/4}%
+1/2\rfloor $ and $u_{\lambda }:=\lfloor \lambda +\sqrt{\lambda
+1/4}+1/2\rfloor $. Moreover,
\begin{eqnarray*}
\bigl\Vert \Delta ^{2}f_{\mathit{Po}(\lambda )}\bigr\Vert _{1}\leq 4\bigl\Vert
\Delta f_{\mathit{Po}(\lambda )}\bigr\Vert _{\infty }  \quad\mbox{and}\quad \bigl\Vert
\Delta ^{2}f_{\mathit{Po}(\lambda )}\bigr\Vert _{1}\sim \frac{4}{\lambda \sqrt{2\curpi \mathrm{e}
}}   \qquad\mbox{as }\lambda \rightarrow \infty .
\end{eqnarray*}
\end{lemma}

\begin{pf}
For convenience, we set $k_{\lambda }:=\lfloor \rho _{1}\rfloor$
and $u_{\lambda }:=\lfloor \rho _{2}\rfloor$, where $\rho
_{1}:=\lambda -\sqrt{\lambda +1/4}+1/2$ and $\rho _{2}:=\lambda +\sqrt{%
\lambda +1/4}+1/2$, and $g(z):=\Delta f_{\mathit{Po}(\lambda )}(z)$. In the proof of
Lemma \ref{prop1}, we have seen that $0=g(-1)\leq g(0)\leq \cdots\leq
g(k_{\lambda })$, while $g(k_{\lambda })\geq g(k_{\lambda }+1)\geq \cdots \geq
g(u_{\lambda })$ and $g(u_{\lambda })\leq g(u_{\lambda }+1)\leq \cdots$. We then have
\begin{eqnarray*}
\bigl\Vert \Delta ^{2}f_{\mathit{Po}(\lambda )}\bigr\Vert _{1}
&=&\sum_{z=0}^{+\infty }|\Delta g(z)|=\sum_{z=0}^{k_{\lambda }}\Delta
g(z)-\sum_{z=k_{\lambda }+1}^{u_{\lambda }}\Delta g(z)+\sum_{z=u_{\lambda
}+1}^{+\infty }\Delta g(z)
\\
&=&\bigl( g( k_{\lambda }) -g(-1)\bigr) -\bigl( g(
u_{\lambda }) -g( k_{\lambda }) \bigr) +\bigl( 0-g(
u_{\lambda }) \bigr)
\\
&=&2\bigl( g( k_{\lambda }) -g( u_{\lambda }) \bigr) .
\end{eqnarray*}
From the proof of Lemma \ref{prop1}, we also get that
\begin{eqnarray*}
g( k_{\lambda }) =\Delta f_{\mathit{Po}(\lambda )}(k_{\lambda
})=\max_{z\in \mathbb{Z}_{+}}\Delta f_{\mathit{Po}(\lambda )}(z)=\bigl\Vert \Delta
f_{\mathit{Po}(\lambda )}\bigr\Vert _{\infty };\qquad g( u_{\lambda })
=\min_{z\in \mathbb{Z}_{+}}\Delta f_{\mathit{Po}(\lambda )}(z)<0
\end{eqnarray*}
and $g( k_{\lambda }) >-g( u_{\lambda }) $. Therefore,
we obtain that $\Vert \Delta ^{2}f_{\mathit{Po}(\lambda )}\Vert _{1}\leq
4g( k_{\lambda }) =4\Vert \Delta f_{\mathit{Po}(\lambda )}\Vert
_{\infty }$. The last asymptotic result follows from the fact that $\Delta
f_{\mathit{Po}(\lambda )}(u_{\lambda })\sim -\lambda ^{-1}(2\curpi \mathrm{e})^{-1/2}$, which can
be proven in exactly the same way as $\Delta f_{\mathit{Po}(\lambda )}(k_{\lambda
})\sim \lambda ^{-1}(2\curpi \mathrm{e})^{-1/2}$ was proven in Lemma \ref{prop1}.
\end{pf}

A crude but simple upper bound is $\Vert \Delta ^{2}f_{\mathit{Po}(\lambda
)}\Vert _{1}\leq 4\frac{1-\mathrm{e}^{-3\lambda }}{3\lambda }\leq 4(1\wedge
\frac{1}{3\lambda })$ for all $\lambda >0$, whereas $\Vert \Delta
f_{\mathit{Po}(\lambda )}\Vert _{\infty }\leq 1/(3\lambda )$ for $\lambda \geq
2 $.

\begin{remark}\label{r1} For distributions other than Poisson, it is not always
easy to derive an analytic expression for $\Vert \Delta f\Vert
_{\infty }$ or $\Vert \Delta ^{2}f\Vert _{1}$. Nevertheless, it
is always feasible to compute the numeric value of these norms employing
numerical or symbolic mathematics software packages (for example, Mathematica, Maple
or MATLAB).
\end{remark}

An approximate expression for these norms can be easily derived if we assume
that the distribution corresponding to the p.d.f. $f$ can be approximated by
a normal distribution $N(\mu ,\sigma ^{2})$, for example, due to CLT. In this case,
we expect that $\Vert \Delta f\Vert _{\infty }$ and $\Vert
\Delta ^{2}f\Vert _{1}$ would be close to $\Vert f_{N(\mu ,\sigma
^{2})}^{( 1) }\Vert _{\infty }$ and $\Vert f_{N(\mu ,\sigma
^{2})}^{( 2) }\Vert _{1}$, respectively, where $f_{N(\mu ,\sigma
^{2})}^{( k) }$ denotes the $k$th order derivative of the p.d.f. of $N(\mu ,\sigma ^{2})$.
It is not difficult to verify that, for the normal
distribution, we have
\begin{eqnarray*}
\bigl\Vert f_{N(\mu ,\sigma ^{2})}^{( 1) }\bigr\Vert_{\infty } &=&\sup_{x\in
\mathbb{R}}\bigl|f_{N(\mu ,\sigma ^{2})}^{( 1) }(x)\bigr|=\frac{1}{\sigma
^{2}\sqrt{2\curpi \mathrm{e}}},
\\
\bigl\Vert f_{N(\mu ,\sigma ^{2})}^{( 2) }\bigr\Vert_{1} &=&\int_{-\infty
}^{+\infty }\bigl|f_{N(\mu ,\sigma ^{2})}^{( 2) }(x)\bigr|\,\mathrm{d}x=4\bigl\Vert f_{N(\mu
,\sigma ^{2})}^{( 1) }\bigr\Vert_{\infty }.
\end{eqnarray*}
Hence, for distributions similar to the normal with variance $\sigma ^{2}$, we
expect $\Vert \Delta ^{2}f\Vert _{1}$ to be nearly equal to $%
4\sigma ^{-2}(2\curpi \mathrm{e})^{-1/2}$. This approximation works for the Poisson
distribution (as seen in Proposition \ref{np} above) since, for large $%
\lambda $, it is close to a normal distribution with $\sigma ^{2}=\lambda $.
According to the above, concerning the compound Poisson distribution, if $%
\mathit{CP}(\lambda ,F)\approx N(\lambda E(W),\lambda E(W^{2}))$ (with $W\sim F$)
then we can expect that, for large $\lambda $,
%e4 ###
\begin{equation}\label{eq25}
\bigl\Vert \Delta ^{2}f_{\mathit{CP}(\lambda ,F)}\bigr\Vert _{1}\approx \frac{4}{%
\lambda E(W^{2})\sqrt{2\curpi \mathrm{e}}}.
\end{equation}

It is worth stressing that (\ref{eq25}) is valid provided  the
compounding distribution $F$ is such that $\mathit{CP}(\lambda ,F)$ is approximately
normal. There exist counterexamples showing that (\ref{eq25}) is not always
valid; see Example $1.3$ or $1.4$ of Barbour and Utev (\citeyear{1999Barbour}).
Specifically, the $\mathit{CP}(\lambda ,F)$ described there cannot be approximated by
a normal distribution and, moreover, it can be verified that the
corresponding $\Vert \Delta ^{2}f_{\mathit{CP}(\lambda ,F)}\Vert _{1}$
does not decrease as $\lambda $ increases. Note that Barbour and Utev (\citeyear{1999Barbour})
use these counterexamples to show that, even for independent $X_{i}$'s (with
$p_{i}=P(X_{i}\neq 0)$, $\lambda =\Sigma p_{i})$, we cannot always prove
that $d_{\mathrm{TV}}(\mathcal{L}(\Sigma X_{i}),\mathit{CP}(\lambda ,F))=\mathrm{O}(\lambda ^{-1}\Sigma
p_{i}^{2})$ and sometimes (depending on $F$) the order $\mathrm{O}(\Sigma p_{i}^{2})$
is optimal. Theorem \ref{iCP} below implies that this $d_{\mathrm{TV}}$ is of order
$\mathrm{O}(\lambda ^{-1}\Sigma p_{i}^{2})$ whenever $F$ is such that $\Vert
\Delta ^{2}f_{\mathit{CP}(\lambda ,F)}\Vert _{1}=\mathrm{O}(\lambda ^{-1})$ (see also
Remark \ref{r2} in Section~\ref{sec3}).\looseness=1

%s2.2 ###
\subsection{Coupling techniques}\label{sec2.2}

A coupling of two random vectors $\mathbf{X},\mathbf{Y}\in\mathbb{R}^{k}$
(to be more exact, of their distributions $\mathcal{L}\mathbf{X}$, $\mathcal{%
L}\mathbf{Y}$) is considered to be any random vector $(\mathbf{X}^{\prime },%
\mathbf{Y}^{\prime })$ defined over a probability space $(\Omega ,\mathfrak{F%
},P)$ and taking values in a measurable space $(\mathbb{R}^{2k},\mathcal{B}(%
\mathbb{R}^{2k}))$ with the same marginal distributions as $\mathbf{X},%
\mathbf{Y}$, that is, $\mathcal{L}\mathbf{X}=\mathcal{L}\mathbf{X}^{\prime } $%
and $\mathcal{L}\mathbf{Y}=\mathcal{L}\mathbf{Y}^{\prime }$. Loosely
speaking, a coupling of $\mathbf{X},\mathbf{Y}$ is any
``definition'' of $\mathbf{X},\mathbf{Y}$ in the same
probability space. This definition of coupling can be generalized for $n$
random vectors in an obvious way. A well-known result concerning the $d_{\mathrm{TV}}$
is the so-called (\textit{basic}) \textit{coupling inequality},
\begin{eqnarray*}
d_{\mathrm{TV}}(\mathcal{L}\mathbf{X},\mathcal{L}\mathbf{Y})\leq P(\mathbf{X}^{\prime
}\neq \mathbf{Y}^{\prime }),
\end{eqnarray*}
which is valid for any coupling $(\mathbf{X}^{\prime },\mathbf{Y}^{\prime })$
of two random vectors $\mathbf{X},\mathbf{Y}$. It can be proven that we can
always construct a coupling $(\mathbf{X}^{\prime },\mathbf{Y}^{\prime })$ of
$( \mathbf{X},\mathbf{Y}) $ such that $d_{\mathrm{TV}}(\mathcal{L}\mathbf{X%
},\mathcal{L}\mathbf{Y})=P(\mathbf{X}^{\prime }\neq \mathbf{Y}^{\prime })$
(for example, see Lindvall (\citeyear{1992Lindvall}), page 18). Such a coupling is called a \textit{%
maximal coupling} or $\gamma $-\textit{coupling} of $\mathbf{X},\mathbf{Y}$.
All of the above could be expressed equivalently for probability measures as
follows: if $P_{1},P_{2}$ are two probability measures on $(\mathbb{R}^{k},%
\mathcal{B}(\mathbb{R}^{k}))$, then any probability measure $\hat{P}$ on $(%
\mathbb{R}^{2k},\mathcal{B}(\mathbb{R}^{2k}))$ with $\hat{P}(A\times
\mathbb{R}^{k})=P_{1}(A)$, $\hat{P}(\mathbb{R}^{k}\times A)=P_{2}(A)$ for
every $A\in \mathcal{B}(\mathbb{R}^{k})$ is called a \textit{coupling} of $%
P_{1},P_{2}$. Moreover, it can be proven that there exists a coupling $\hat{P%
}_{\gamma }$ of $P_{1},P_{2}$, called a \textit{maximal coupling} or $\gamma $\emph{-}%
\textit{coupling}, such that
%e5 ###
\begin{equation} \label{couple}
d_{\mathrm{TV}}(P_{1},P_{2})=1-\hat{P}_{\gamma }\bigl(\{(\mathbf{x},\mathbf{x}),\mathbf{x}%
\in \mathbb{R}^{k}\}\bigr).
\end{equation}
Obviously, all of the above can be adapted in the obvious way for random
vectors taking values in $\mathbb{Z}^{k}$ and to multivariate distributions
over the probability space $(\mathbb{Z}^{k},2^{\mathbb{Z}^{k}})$.

The following lemmas will play a crucial role for the establishment of our
main results. The first inequality of the following lemma is Corollary $%
4 $ in Boutsikas (\citeyear{2006Boutsikas}). The second inequality of the following lemma is a
direct application of Lemma $3$ in Boutsikas (\citeyear{2006Boutsikas}) with $(\Xi _{1}^{\prime
},\Xi _{2}^{\prime },\Psi _{1}^{\prime },\Psi _{2}^{\prime })=(\mathbf{Z}+%
\mathbf{X},\mathbf{Z}+\mathbf{Y},\mathbf{W}+\mathbf{X},\mathbf{W}+\mathbf{Y)}
$.

\begin{lemma}
\label{lemma3} For any random vectors $\mathbf{X},\mathbf{Y}\in \mathbb{R}%
^{k}$ and $\mathbf{Z},\mathbf{W}\in \mathbb{R}^{r}$ defined on the same
probability space, we have that
\begin{longlist}
\item[(a)] $\vert d_{\mathrm{TV}}(\mathcal{L}(\mathbf{Z},\mathbf{X}),\mathcal{L}(%
\mathbf{Z},\mathbf{Y}))-d_{\mathrm{TV}}(\mathcal{L}(\mathbf{W},\mathbf{X}),\mathcal{L}%
(\mathbf{W},\mathbf{Y}))\vert \leq 2P( \mathbf{X}\neq \mathbf{Y,Z}%
\neq \mathbf{W});$

\item[(b)] $\vert d_{\mathrm{TV}}(\mathcal{L}(\mathbf{Z}+\mathbf{X)},\mathcal{L}(%
\mathbf{Z}+\mathbf{Y)})-d_{\mathrm{TV}}(\mathcal{L}(\mathbf{W}+\mathbf{X)},\mathcal{L}%
(\mathbf{W}+\mathbf{Y)})\vert \leq 2P(\mathbf{X}\neq \mathbf{Y},%
\mathbf{Z}\neq \mathbf{W}).$
\end{longlist}
\end{lemma}

The next inequality follows from the above lemma. It is remarkable that
almost the same inequality can be found in Rachev (\citeyear{1991Rachev}), page 274, and has
been applied to derive Berry--Esseen-type results. We present an entirely
different proof using maximal couplings.

\begin{lemma}
\label{lemma2} If the random vectors $\mathbf{X},\mathbf{Y}\in \mathbb{R}%
^{k} $are independent of $\mathbf{Z},\mathbf{W}\in \mathbb{R}^{r}$, then%
\begin{eqnarray*}
&& \bigl\vert d_{\mathrm{TV}}\bigl(\mathcal{L}(\mathbf{Z}+\mathbf{X)},\mathcal{%
L}(\mathbf{Z}+\mathbf{Y)}\bigr)-d_{\mathrm{TV}}\bigl(\mathcal{L}(\mathbf{W}+\mathbf{X)},%
\mathcal{L}(\mathbf{W}+\mathbf{Y)}\bigr)\bigr\vert
\\
&&\quad \leq 2d_{\mathrm{TV}}(\mathcal{L}\mathbf{X},\mathcal{L}\mathbf{Y})d_{\mathrm{TV}}(\mathcal{L}%
\mathbf{Z},\mathcal{L}\mathbf{W}).
\end{eqnarray*}
\end{lemma}

\begin{pf}
Let $(\mathbf{X}^{\ast },\mathbf{Y}^{\ast })$ be a maximal coupling of $%
\mathcal{L}\mathbf{X},\mathcal{L}\mathbf{Y}$ and let $(\mathbf{Z}^{\ast },%
\mathbf{W}^{\ast })$ be a maximal coupling of $\mathcal{L}\mathbf{Z}$, $%
\mathcal{L}\mathbf{W}$. Next, let $((\mathbf{X}^{\prime },\mathbf{Y}^{\prime
})$, $(\mathbf{Z}^{\prime },\mathbf{W}^{\prime }))$ be an independent
coupling of $\mathcal{L}(\mathbf{X}^{\ast },\mathbf{Y}^{\ast })$, $\mathcal{L%
}(\mathbf{Z}^{\ast },\mathbf{W}^{\ast })$ (that is, $(\mathbf{X}^{\prime },%
\mathbf{Y}^{\prime })$ is independent of $(\mathbf{Z}^{\prime },\mathbf{W}%
^{\prime })$ and $\mathcal{L}(\mathbf{X}^{\prime },\mathbf{Y}^{\prime })=%
\mathcal{L}(\mathbf{X}^{\ast },\mathbf{Y}^{\ast })$, $\mathcal{L}(\mathbf{Z}%
^{\prime },\mathbf{W}^{\prime })=\mathcal{L}(\mathbf{Z}^{\ast },\mathbf{W}%
^{\ast })$). Applying Lemma \ref{lemma3}, we get
\begin{eqnarray*}
&&\bigl\vert d_{\mathrm{TV}}\bigl(\mathcal{L}(\mathbf{Z}^{\prime }+\mathbf{X}%
^{\prime }\mathbf{)},\mathcal{L}(\mathbf{Z}^{\prime }+\mathbf{Y}^{\prime }%
\mathbf{)}\bigr)-d_{\mathrm{TV}}\bigl(\mathcal{L}(\mathbf{W}^{\prime }+\mathbf{X}^{\prime }%
\mathbf{)},\mathcal{L}(\mathbf{W}^{\prime }+\mathbf{Y}^{\prime }\mathbf{)}%
\bigr)\bigr\vert
\\
&&\quad \leq 2P(\mathbf{X}^{\prime }\neq \mathbf{Y}^{\prime },\mathbf{Z}^{\prime
}\neq \mathbf{W}^{\prime })=2P(\mathbf{X}^{\prime }\neq \mathbf{Y}^{\prime
})P(\mathbf{Z}^{\prime }\neq \mathbf{W}^{\prime })
\\
&&\quad =2P(\mathbf{X}^{\ast }\neq \mathbf{Y}^{\ast })P(\mathbf{Z}^{\ast }\neq
\mathbf{W}^{\ast })=2d_{\mathrm{TV}}(\mathcal{L}\mathbf{X},\mathcal{L}\mathbf{Y}%
)d_{\mathrm{TV}}(\mathcal{L}\mathbf{Z},\mathcal{L}\mathbf{W}).
\end{eqnarray*}
The obvious fact that $\mathcal{L}(\mathbf{Z}^{\prime }+\mathbf{X}^{\prime }%
)=\mathcal{L}(\mathbf{Z}+\mathbf{X)}$, $\mathcal{L}(\mathbf{Z}%
^{\prime }+\mathbf{Y}^{\prime })=\mathcal{L}(\mathbf{Z}+\mathbf{Y)}$%
, $\mathcal{L}(\mathbf{W}^{\prime }+\mathbf{X}^{\prime })=\mathcal{L%
}(\mathbf{W}+\mathbf{X)}$ and $\mathcal{L}(\mathbf{W}^{\prime }+\mathbf{Y}%
^{\prime })=\mathcal{L}(\mathbf{W}+\mathbf{Y)}$ completes the proof.
\end{pf}

A direct application of the previous result leads to the following
inequality which is valid for any random variables $X,Y\in \mathbb{R} $
independent of another random variable $W\in \mathbb{R}$. Specifically, if
we simply set $Z=0$ in Lemma \ref{lemma2} and exploit the fact that $d_{\mathrm{TV}}(%
\mathcal{L}0,\mathcal{L}W)=P(W\neq 0)$, we derive
\begin{eqnarray*}
d_{\mathrm{TV}}(\mathcal{L}X,\mathcal{L}Y)\leq 2d_{\mathrm{TV}}\bigl(\mathcal{L}X,\mathcal{L}%
Y\bigr)P(W\neq 0)+d_{\mathrm{TV}}\bigl(\mathcal{L(}X+W),\mathcal{L}(Y+W)\bigr)
\end{eqnarray*}
which, for $P(W\neq 0)<1/2$, implies that
%e6 ###
\begin{equation}
d_{\mathrm{TV}}(\mathcal{L}X,\mathcal{L}Y)\leq \bigl( 1-2P(W\neq 0)\bigr)
^{-1} d_{\mathrm{TV}}\bigl(\mathcal{L}(X+W),\mathcal{L}(Y+W)\bigr).  \label{eq17}
\end{equation}

The next lemma can be considered as a coupling inequality concerning \textbf{$%
\zeta $}$_{2}$, analogous to Lemma \ref{lemma3}.

\begin{lemma}
\label{th5}If $X,Y,Z,W$ are real-valued, non-negative random variables
defined on the same probability space with finite second moments and $%
E(X)=E(Y)$, then
%e7 ###
\begin{equation}\label{eq44}
\bigl\vert \mathbf{\zeta }_{2}\bigl(\mathcal{L}(X+Z),\mathcal{L}(Y+Z)\bigr)-\mathbf{%
\zeta }_{2}\bigl(\mathcal{L}(X+W),\mathcal{L}(Y+W)\bigr)\bigl\vert \leq E|(X-Y)(Z-W)|.
\end{equation}
\end{lemma}

\begin{pf}
The distances $\mathbf{\zeta }_{2}$ appearing in (\ref{eq44}) are well
defined since the random variables $X+Z$, $Y+Z$, $X+W$ and $Y+W$ have finite
second moments due to Minkowski's inequality and $E(X+Z)=E(Y+Z)$, $%
E(X+W)=E(Y+W)$. Set $\mathbf{1}_{[a\leq b]}:=1$ if $a\leq b$ and $\mathbf{1}%
_{[a\leq b]}:=0$ otherwise. As usual, $F_{V}$ denotes the c.d.f. of a random
variable $V$. Recall that, for $X,Y\in \mathbb{R}_{+}$ with $E(X)=E(Y)$, we
have
\begin{eqnarray*}
\mathbf{\zeta }_{2}(\mathcal{L}X,\mathcal{L}Y)=\int_{0}^{\infty }\vert
E(X-s)_{+}-E(Y-s)_{+}\vert \,\mathrm{d}s=\int_{0}^{\infty }\biggl\vert
\int_{s}^{\infty }\bigl( F_{X}(x)-F_{Y}(x)\bigr) \,\mathrm{d}x\biggr\vert\, \mathrm{d}s.
\end{eqnarray*}
Denoting by $d$ the absolute difference in the left-hand side of (\ref{eq44}), we
have
\begin{eqnarray*}
d &=&\biggl\vert \int_{0}^{\infty }\biggl\vert \int_{s}^{\infty }\bigl(
F_{X+Z}(x)-F_{Y+Z}(x)\bigr)\, \mathrm{d}x\biggr\vert\,\mathrm{d}s-\int_{0}^{\infty }\biggl\vert
\int_{s}^{\infty }\bigl( F_{X+W}(x)-F_{Y+W}(x)\bigr) \,\mathrm{d}x\biggr\vert\,
\mathrm{d}s\biggr\vert
\\
&\leq &\int_{0}^{\infty }\biggl\vert \biggl\vert\int_{s}^{\infty }\bigl(
F_{X+Z}(x)-F_{Y+Z}(x)\bigr) \,\mathrm{d}x\biggr\vert -\biggl\vert \int_{s}^{\infty
}\bigl( F_{X+W}(x)-F_{Y+W}(x)\bigr) \,\mathrm{d}x \biggr\vert\biggr\vert \,\mathrm{d}s.
\end{eqnarray*}
Using the inequality $\vert \vert a\vert -\vert
b\vert \vert \leq \vert a-b\vert , a,b\in \mathbb{%
\mathbb{R}}$, we get
\begin{eqnarray*}
d &\leq &\int_{0}^{\infty }\biggl\vert \int_{s}^{\infty }\bigl(
F_{X+Z}(x)-F_{Y+Z}(x)\bigr)\, \mathrm{d}x-\int_{s}^{\infty }\bigl(
F_{X+W}(x)-F_{Y+W}(x)\bigr)\, \mathrm{d}x\biggr\vert\, \mathrm{d}s
\\
&=&\int_{0}^{\infty }\vert E( C_{s}) \vert \,\mathrm{d}s\leq
E\biggl( \int_{0}^{\infty }\vert C_{s}\vert \,\mathrm{d}s\biggr),
\end{eqnarray*}
where
\begin{eqnarray*}
C_{s}=\int_{s}^{\infty }\bigl(\mathbf{1}_{[X+Z\leq x]}-\mathbf{1}_{[Y+Z\leq
x]}\bigr)\,\mathrm{d}x-\int_{s}^{\infty }\bigl(\mathbf{1}_{[X+W\leq x]}-\mathbf{1}_{[Y+W\leq
x]}\bigr)\,\mathrm{d}x.
\end{eqnarray*}
Now, if $Z\geq W$, it can be verified that $C_{s}\geq 0$ for all $s>0$ and,
therefore,%
\begin{eqnarray*}
\int_{0}^{\infty }\vert C_{s}\vert\, \mathrm{d}s &=&\int_{0}^{\infty }x\bigl(%
\mathbf{1}_{[X+Z\leq x]}-\mathbf{1}_{[Y+Z\leq x]}\bigr)\,\mathrm{d}x-\int_{0}^{\infty }x\bigl(%
\mathbf{1}_{[X+W\leq x]}-\mathbf{1}_{[Y+W\leq x]}\bigr)\,\mathrm{d}x
\\
&=&\frac{1}{2}\bigl( \vert (X+Z)^{2}-(Y+Z)^{2}\vert -\vert
(X+W)^{2}-(Y+W)^{2}\vert \bigr)
\\
&=&\frac{1}{2}\bigl( \vert (X-Y)(X+Y+2Z)\vert -\vert
(X-Y)(X+Y+2W)\vert \bigr)
\\
&=&|X-Y|(Z-W).
\end{eqnarray*}
On the other hand, if $Z\leq W$, then $C_{s}\leq 0$ for all $s>0$ and we
similarly derive that $\int_{0}^{\infty }\vert C_{s}\vert
\,\mathrm{d}s=|X-Y|(W-Z)$. Hence, $\int_{0}^{\infty }\vert C_{s}\vert
\,\mathrm{d}s=|(X-Y)(W-Z)|$ and the proof is completed.
\end{pf}

A direct corollary of Lemma \ref{th5} is the following result which will be
proven useful when dealing with $k$-dependent sequences of random variables.

\begin{corollary}
\label{cor4}If the random variables $X_{1},X_{2},\ldots,X_{i}\in \mathbb{R}_{+}$
are $k$-dependent with $E(X_{i}^{2})<\infty $ and $l\leq i-k+1$, then
\begin{eqnarray*}
\mathbf{\zeta }_{2}\Biggl( \mathcal{L}\sum_{j=l}^{i}X_{j},\mathcal{L}\Biggl(
\sum_{j=l}^{i-1}X_{j}+X_{i}^{\perp }\Biggr) \Biggr) \leq
\sum_{j=i-k+1}^{i-1}\bigl( E(X_{i}X_{j})+E(X_{i})E(X_{j})\bigr),
\end{eqnarray*}
where $X_{i}^{\perp }$ is a random variable independent of all $X_{j}$, $%
j=1,2,\ldots,i$, with $\mathcal{L}X_{i}=\mathcal{L}X_{i}^{\perp }$.
\end{corollary}

\begin{pf}
Set $X_{a,b}:=\sum_{j=a}^{b}X_{j}$. Applying Lemma \ref{th5} with $%
X=X_{i},Z=X_{l,i-1},Y=X_{i}^{\perp },W=X_{l,i-k}$, we obtain
\begin{eqnarray*}
&&\bigl\vert \mathbf{\zeta }_{2}\bigl(\mathcal{L}X_{l,i},\mathcal{L}%
(X_{l,i-1}+X_{i}^{\perp })\bigr)-\mathbf{\zeta }_{2}\bigl(\mathcal{L}(X_{l,i-k}+X_{i}),%
\mathcal{L}(X_{l,i-k}+X_{i}^{\perp })\bigr)\bigr\vert
\\
&&\quad \leq E|(X_{i}-X_{i}^{\perp })(X_{l,i-1}-X_{l,i-k})|=E|(X_{i}-X_{i}^{\perp
})(X_{i-k+1,i-1})|
\\
&&\quad \leq \sum_{j=i-k+1}^{i-1}\bigl( E(X_{i}X_{j})+E(X_{i})E(X_{j})\bigr).
\end{eqnarray*}
Since $X_{l,i-k}$ and $X_{i}$ are independent, $X_{l,i-k}$ and $%
X_{i}^{\perp }$ are independent, and $\mathcal{L}X_{i}=\mathcal{L}X_{i}^{\perp }$, we
conclude that $\mathcal{L}(X_{l,i-k}+X_{i})=\mathcal{L}(X_{l,i-k}+X_{i}^{%
\perp })$ and hence we  obtain the desired inequality.
\end{pf}

As will be seen in the next section, Lemmas \ref{Btv} and \ref{lemma2}
 are sufficient for proving compound Poisson approximation results
for sums of \textit{independent} random variables incorporating a smoothness
factor. In the case of sums of dependent random variables, though, the
following, additional, lemma is needed. The question addressed here is the
following: given a random variable $X$ and a random vector $\mathbf{Z}$, can
we construct (on the same probability space as $X,\mathbf{Z)}$ another
random variable $Y$ with a given p.d.f. $f$ such that $Y$ is independent of $%
\mathbf{Z}$ and $(X,\mathbf{Z}),(Y,\mathbf{Z})$ are maximally coupled? In
this situation, we could loosely say that we wish to construct a random
variable $Y$ (with a given distribution) that  resembles $X$ as far as
possible, \textit{while remaining independent of~$\mathbf{Z}$}. Again, it
suffices to restrict our analysis to the discrete case.

\begin{lemma}
\label{lemma6} Let $X\in \mathbb{Z},\mathbf{Z}\in \mathbb{Z}^{k}$ be a
random variable and a random vector, respectively (defined on the same
probability space) and let $f\dvtx\mathbb{Z\rightarrow R}_{+}$ be some given
discrete p.d.f. Denote by $U$ a random variable independent of $X,\mathbf{Z}
$ that follows the uniform distribution on $(0,1)$. Then,
\begin{longlist}
\item[(a)] there exists a function $g\dvtx \mathbb{R}^{2+k}\rightarrow \mathbb{R}$ such
that the random variable $Y=g(U,X,\mathbf{Z})$ has p.d.f. $f$, $Y$ is
independent of $\mathbf{Z}$ and
\begin{eqnarray*}
d_{\mathrm{TV}}(\mathcal{L}(X,\mathbf{Z}),\mathcal{L}(Y,\mathbf{Z}))=P\bigl((X,\mathbf{Z}%
)\neq (Y,\mathbf{Z})\bigr)=P( X\neq Y),
\end{eqnarray*}
in other words, $(X,\mathbf{Z}),(Y,\mathbf{Z})$ are maximally coupled;

\item[(b)] there exists a function $g^{\prime }\dvtx\mathbb{R}^{2}\rightarrow \mathbb{R}
$ such that the random variable $Y^{\prime }=g^{\prime }(U,X)$ has p.d.f.~$f$
and $(X,Y^{\prime })$ are maximally coupled, that is, $d_{\mathrm{TV}}(\mathcal{L}X,%
\mathcal{L}Y^{\prime })=P(X\neq Y^{\prime })$.
\end{longlist}
\end{lemma}

\begin{pf}
(a) Here, we develop a constructive proof. Denote by $(\Omega ,\mathcal{A}%
,P)$ the probability space on which $X,\mathbf{Z,}U$ are defined and let $%
f_{X|\mathbf{Z}}(\cdot |\mathbf{z})=f_{X,\mathbf{Z}}(\cdot ,\mathbf{z})/f_{%
\mathbf{Z}}(\mathbf{z})$ be the conditional p.d.f. of $X$ given $\mathbf{Z}=%
\mathbf{z}$. Consider the probability measures $P_{1}^{\mathbf{z}},P_{2}$ on
the measurable space $(\mathbb{Z},2^{\mathbb{Z}})$ generated by $f_{X|%
\mathbf{Z}}(\cdot |\mathbf{z})$ and $f$, respectively. According to (\ref%
{couple}), there exists a maximal coupling of $P_{1}^{\mathbf{z}},P_{2}$.
Denote by $h_{\mathbf{z}}\dvtx\mathbb{Z}^{2}\rightarrow \mathbb{R}_{+}$ the
joint p.d.f. corresponding to this maximal coupling. It follows that $%
\sum_{x\in \mathbb{Z}}h_{\mathbf{z}}(x\mathbf{,}y)=f(y)$, $\sum_{y\in
\mathbb{Z}}h_{\mathbf{z}}(x\mathbf{,}y)=f_{X|\mathbf{Z}}(x|\mathbf{z})$ and
\begin{eqnarray*}
d_{\mathrm{TV}}\bigl(\mathcal{L}(X|\mathbf{Z}=\mathbf{z)},P_{2}\bigr)=d_{\mathrm{TV}}(P_{1}^{\mathbf{z}%
},P_{2})=1-\sum_{x\in \mathbb{Z}}h_{\mathbf{z}}(x\mathbf{,}x).
\end{eqnarray*}
We now construct $Y$ as follows. For every $x\in \mathbb{Z},\mathbf{z}\in
\mathbb{Z}^{k}$, consider the c.d.f.
\begin{eqnarray*}
H_{x,\mathbf{z}}(y):=\sum_{i\leq \lfloor y\rfloor }\frac{h_{%
\mathbf{z}}(x,i)}{f_{X|\mathbf{Z}}(x|\mathbf{z})},\qquad y\in \mathbb{R},
\end{eqnarray*}%
and set $Y(\omega ):=H_{X(\omega ),\mathbf{Z}(\omega )}^{-1}(U(\omega ))$, $%
\omega \in \Omega $, where $H_{x,\mathbf{z}}^{-1}(y)$ denotes the
generalized inverse of $H_{x,\mathbf{z}}(y)$, that is, $H_{x,\mathbf{z}%
}^{-1}(y)=\inf \{ w\dvtx H_{x,\mathbf{z}}(w)\geq y\} $. The function $%
f_{X,Y,\mathbf{Z}}(x,y,\mathbf{z}):=h_{\mathbf{z}}(x\mathbf{,}y)f_{\mathbf{Z}%
}(\mathbf{z})$ is a multivariate discrete p.d.f. and it can be verified that
$Y$ and $(X,Y,\mathbf{Z)}$ have p.d.f. $f$ and $f_{X,Y,\mathbf{Z}}$,
respectively. Indeed,
\begin{eqnarray*}
P( X=x,Y\leq y,\mathbf{Z}=\mathbf{z}) &=&P\bigl(X=x,H_{x,\mathbf{z}%
}^{-1}(U)\leq y,\mathbf{Z}=\mathbf{z}\bigr)
\\
&=&H_{x,\mathbf{z}}(y)P(X=x,\mathbf{Z}=\mathbf{z})
\end{eqnarray*}
and thus, for all $x,y,z$,
\begin{eqnarray*}
P( X=x,Y=y,\mathbf{Z}=\mathbf{z}) &=&\frac{h_{\mathbf{z}}(x%
\mathbf{,}y)}{f_{X|\mathbf{Z}}(x|\mathbf{z})}P(X=x,\mathbf{Z}=\mathbf{z})
\\
&=&h_{\mathbf{z}}(x\mathbf{,}y)f_{\mathbf{Z}}(\mathbf{z})=f_{X,Y,\mathbf{Z}
}(x,y,\mathbf{z}).
\end{eqnarray*}
Also, note that, for all $x,\mathbf{z}$,
\begin{eqnarray*}
P( Y=y,\mathbf{Z}=\mathbf{z}) =\sum_{x\in \mathbb{Z}}f_{X,Y,%
\mathbf{Z}}(x,y,\mathbf{z})=\sum_{x\in \mathbb{Z}}h_{\mathbf{z}}(x,
y)f_{\mathbf{Z}}(\mathbf{z})=f(y)f_{\mathbf{Z}}(\mathbf{z}),
\end{eqnarray*}
which implies that $Y$ is independent of  $ \mathbf{Z}$.
Furthermore, we derive that, for all $\mathbf{z}$,
\begin{eqnarray*}
P(X\neq Y|\mathbf{Z}=\mathbf{z})=1-\sum_{x\in \mathbb{Z}}h_{\mathbf{z}}(x
,x)=d_{\mathrm{TV}}\bigl(\mathcal{L}(X|\mathbf{Z}=\mathbf{z)},\mathcal{L}Y\bigr)
\end{eqnarray*}
and, therefore,
\begin{eqnarray*}
P( X\neq Y) &=&\sum_{\mathbf{z\in }\mathbb{Z}^{k}}P(X\neq Y|%
\mathbf{Z}=\mathbf{z})f_{\mathbf{Z}}(\mathbf{z})
\\
&=&\sum_{\mathbf{z\in }\mathbb{Z}^{k}}d_{\mathrm{TV}}\bigl(\mathcal{L}(X|\mathbf{Z}=%
\mathbf{z)},\mathcal{L}Y\bigr)f_{\mathbf{Z}}(\mathbf{z})
\\
&=&\sum_{\mathbf{z\in }\mathbb{Z}^{k}}\frac{1}{2}\sum_{w\in \mathbb{Z}%
}\vert P(X=w|\mathbf{Z}=\mathbf{z})-P(Y=w)\vert f_{\mathbf{Z}}(%
\mathbf{z}) \\
&=&\frac{1}{2}\sum_{\mathbf{z\in }\mathbb{Z}^{k}}\sum_{w\in \mathbb{Z}%
}\vert P(X=w,\mathbf{Z}=\mathbf{z})-P(Y=w)f_{\mathbf{Z}}(\mathbf{z}%
)\vert \\
&=&d_{\mathrm{TV}}(\mathcal{L}(X,\mathbf{Z}),\mathcal{L}(Y,\mathbf{Z})).
\end{eqnarray*}

(b) This  readily follows from part (a) of the lemma by choosing $\mathbf{Z}=%
\mathbf{0}$.
\end{pf}

%s3 ###
\section{Compound Poisson approximation for sums of independent random
variables}\label{sec3}

Let $X_{1},X_{2},\ldots ,X_{n}$ be a sequence of independent random variables
which take values in $\mathbb{Z}_{+}$. We are now ready to exploit the
results of the previous section (specifically Lemmas \ref{Btv} and \ref%
{lemma2}) to derive a simple and, in most cases, sharp upper bound for the
total variation distance between the distribution of the sum $%
\sum_{i=1}^{n}X_{i}$ and an appropriate compound Poisson distribution.
Before we present this bound, we recall that, (see Boutsikas and Vaggelatou
(\citeyear{2002Boutsikas}))
%e8 ###
\begin{equation}\label{eq36}
\mathbf{\zeta }_{2}\Biggl( \mathcal{L}\sum_{i=1}^{n}X_{i},\mathit{CP}\Biggl( \lambda ,%
\frac{1}{\lambda }\sum_{i=1}^{n}p_{i}G_{i}\Biggr) \Biggr) =\frac{1}{2}%
\sum_{i=1}^{n}E(X_{i})^{2},
\end{equation}
with $p_{i}:=P(X_{i} \neq 0)$, $\lambda =\sum_{i=1}^{n}p_{i}$ and $%
G_{i}(x)=P(X_{i}\leq x|X_{i}\neq 0)$. Naturally, the bound of the following
theorem is useful (that is, it tends to $0$) when $p_{i}\approx 0$. Hence, the
condition $p_{i}<\log 2\approx 0.693$ imposed below does not affect the
generality of the result. One could easily modify the upper bound (making it
a little bit more complicated) so as to eliminate this restriction, but this
modification would lead to no practical gain.

\begin{theorem}
\label{iCP}Let $X_{1},X_{2},\ldots ,X_{n}$ be a sequence of independent random
variables (with finite second moments) taking values in $\mathbb{Z}_{+}$ and
$P(X_{i}\neq 0)=:p_{i}<\log 2$ $(\approx 0.693)$. Then,
\begin{eqnarray*}
&&d_{\mathrm{TV}}\Biggl( \mathcal{L}\sum_{i=1}^{n}X_{i},\mathit{CP}(\lambda
,F)\Biggr)
\\
&&\quad \leq \Biggl( \sum_{i=1}^{n}p_{i}^{2}\Biggr) ^{2}+\frac{1}{4}\bigl\Vert
\Delta ^{2}f_{\mathit{CP}(\lambda ,F)}\bigr\Vert _{1}\sum_{i=1}^{n}\frac{E(X_{i})^{2}%
}{1-2( 1-\mathrm{e}^{-p_{i}}) }:=UB_{\mathit{CP}},
\end{eqnarray*}
where $\lambda =\sum_{i=1}^{n}p_{i}$, $F(x)=\sum_{i=1}^{n}\frac{p_{i}}{%
\lambda }G_{i}(x)$ and $G_{i}(x)=P(X_{i}\leq x|X_{i}\neq 0)$, $x\in \mathbb{Z%
}$.
\end{theorem}

\begin{pf}
Let $N_{1},N_{2},\ldots ,N_{n}$ be independent random variables following the
compound Poisson distribution with parameters $%
(p_{1},G_{1}),(p_{2},G_{2}),\ldots ,(p_{n},G_{n}),$ respectively. We apply the
triangle inequality to get the following Lindeberg decomposition of the
distance of interest,
%
%e9 ###
\begin{equation} \label{eq19}
d_{\mathrm{TV}}\Biggl( \mathcal{L}\sum_{i=1}^{n}X_{i},\mathcal{L}\sum_{i=1}^{n}N_{i}%
\Biggr) \leq \sum_{m=1}^{n}d_{\mathrm{TV}}\Biggl( \mathcal{L}\Biggl(
\sum_{i=1}^{m}X_{i}+\sum_{i=m+1}^{n}N_{i}\Biggr) ,\mathcal{L}\Biggl(
\sum_{i=1}^{m-1}X_{i}+\sum_{i=m}^{n}N_{i}\Biggr) \Biggr).
\end{equation}
Furthermore, if we set
\begin{eqnarray*}
\mathbf{X}_{m}:=X_{m}+\sum_{i=m+1}^{n}N_{i},\qquad \mathbf{Y}_{m}:=N_{m}+%
\sum_{i=m+1}^{n}N_{i},\qquad \mathbf{Z}_{m}:=\sum_{i=1}^{m-1}X_{i},\qquad \mathbf{W}%
_{m}:=\sum_{i=1}^{m-1}N_{i},
\end{eqnarray*}
then the random variables $\mathbf{X}_{m},\mathbf{Y}_{m} $are independent
of $\mathbf{Z}_{m},\mathbf{W}_{m}$ and a direct application of Lemma~\ref%
{lemma2} to $\mathbf{X}_{m},\mathbf{Y}_{m},\mathbf{Z}_{m},\mathbf{W}_{m}$
reveals that
%
%e10 ###
\begin{equation}\label{eq20}
d_{\mathrm{TV}}\Biggl( \mathcal{L}\Biggl(
\sum_{i=1}^{m-1}X_{i}+X_{m}+\sum_{i=m+1}^{n}N_{i}\Biggr) ,\mathcal{L}\Biggl(
\sum_{i=1}^{m-1}X_{i}+N_{m}+\sum_{i=m+1}^{n}N_{i}\Biggr) \Biggr) \leq
2a_{m}b_{m}+c_{m},
\end{equation}
where
\begin{eqnarray*}
a_{m} &:=&d_{\mathrm{TV}}\Biggl( \mathcal{L}\Biggl( X_{m}+\sum_{i=m+1}^{n}N_{i}\Biggr) ,%
\mathcal{L}\Biggl( N_{m}+\sum_{i=m+1}^{n}N_{i}\Biggr) \Biggr),
\\
b_{m} &:=&d_{\mathrm{TV}}\Biggl( \mathcal{L}\sum_{i=1}^{m-1}X_{i},\mathcal{L}%
\sum_{i=1}^{m-1}N_{i}\Biggr),
\\
c_{m} &:=&d_{\mathrm{TV}}\Biggl( \mathcal{L}\Biggl(
\sum_{i=1}^{n}N_{i}-N_{m}+X_{m}\Biggr) ,\mathcal{L}\sum_{i=1}^{n}N_{i}%
\Biggr).
\end{eqnarray*}
Next, let $N_{m}^{\perp }$ be a random variable independent of all $%
N_{i},X_{i}$ with $\mathcal{L}N_{m}^{\perp }=\mathcal{L}N_{m}$. Applying
inequality (\ref{eq17}) with $W=N_{m}^{\perp }$, we derive
\begin{eqnarray*}
c_{m}\leq \bigl(1-2(1-\mathrm{e}^{-p_{m}})\bigr)^{-1}d_{\mathrm{TV}}\Biggl( \mathcal{L}\Biggl(
\sum\limits_{i=1}^{n}N_{i}+X_{m}\Biggr) ,\mathcal{L}\Biggl(
\sum\limits_{i=1}^{n}N_{i}+N_{m}^{\perp }\Biggr) \Biggr)
\end{eqnarray*}
since $P(N_{m}^{\perp }\neq 0)=1-\mathrm{e}^{-p_{m}}$. Furthermore, Lemma \ref{Btv}
yields
%
%e11 ###
\begin{equation}\label{eq22}
c_{m}\leq \frac{{1/2}\Vert \Delta ^{2}f_{\Sigma
_{i=1}^{n}N_{i}}\Vert _{1}}{1-2(1-\mathrm{e}^{-p_{m}})}\mathbf{\zeta }_{2}(%
\mathcal{L}X_{m},\mathcal{L}N_{m}^{\perp })=\frac{\Vert \Delta
^{2}f_{\mathit{CP}(\lambda ,F)}\Vert _{1}}{4(1-2(1-\mathrm{e}^{-p_{m}}))}E(X_{m})^{2},
\end{equation}
where we have used (\ref{eq36}) to get that $\mathbf{\zeta }_{2}(\mathcal{L}%
X_{m},\mathcal{L}N_{m}^{\perp })=\mathbf{\zeta }_{2}(\mathcal{L}%
X_{m},\mathit{CP}(p_{m},G_{m}))=\frac{1}{2}E(X_{m})^{2}$. On the other hand, we can
easily bound the quantities $a_{m},b_{m}$ as follows:
%
%e12 ###
\begin{equation} \label{eq23}
a_{m}\leq d_{\mathrm{TV}}(\mathcal{L}X_{m},\mathcal{L}N_{m})\leq p_{m}^{2}
\quad \mbox{and}\quad   b_{m}\leq \sum_{i=1}^{m-1}d_{\mathrm{TV}}(\mathcal{L}X_{i},\mathcal{L}%
N_{i})=\sum_{i=1}^{m-1}p_{i}^{2}.
\end{equation}
Finally, combining (\ref{eq19})--(\ref{eq23}), we get
\begin{eqnarray*}
d_{\mathrm{TV}}\Biggl( \mathcal{L}\sum_{i=1}^{n}X_{i},\mathcal{L}\sum_{i=1}^{n}N_{i}%
\Biggr) &\leq &\sum_{m=1}^{n}( 2a_{m}b_{m}+c_{m})
\\
&\leq &\sum_{m=1}^{n}\Biggl( 2p_{m}^{2}\sum_{i=1}^{m-1}p_{i}^{2}+\frac{%
\Vert \Delta ^{2}f_{\mathit{CP}(\lambda ,F)}\Vert _{1}}{%
4(1-2(1-\mathrm{e}^{-p_{m}}))}E(X_{m})^{2}\Biggr),
\end{eqnarray*}
which readily leads to the desired inequality since
\begin{eqnarray*}
\hspace{6pt}\quad\qquad\qquad\sum_{m=1}^{n}2a_{m}b_{m}&\leq&
\sum_{m=1}^{n}2p_{m}^{2}\sum_{i=1}^{m-1}p_{i}^{2}
\\
&=&\sum_{m=1}^{n}p_{m}^{2}%
\sum_{i=1}^{m-1}p_{i}^{2}+\sum_{m=1}^{n}p_{m}^{2}\sum_{i=m+1}^{n}p_{i}^{2}%
\leq \Biggl( \sum_{m=1}^{n}p_{m}^{2}\Biggr) ^{2}.\quad\hspace{26pt}\qed
\end{eqnarray*}
\noqed\end{pf}

A straightforward corollary of the above theorem arises when we consider
independent Bernoulli random variables. In this case, the distribution of
the sum of the binary sequence $X_{1},X_{2},\ldots ,X_{n}$ is also known as a
Poisson binomial or generalized binomial distribution and the approximating
compound Poisson distribution naturally reduces to an ordinary
Poisson
distribution.

\begin{corollary}
\label{iP}Let $X_{1},X_{2},\ldots,X_{n}$ be a sequence of independent Bernoulli
random variables with $P(X_{i}=1)=p_{i}<\log 2,i=1,2,\ldots,n$. Then,
\begin{eqnarray*}
d_{\mathrm{TV}}\Biggl( \sum_{i=1}^{n}X_{i},\mathit{Po}(\lambda )\Biggr)
 \leq \Biggl( \sum_{i=1}^{n}p_{i}^{2}\Biggr) ^{2}+\frac{1}{4}\bigl\Vert
\Delta ^{2}f_{\mathit{Po}(\lambda )}\bigr\Vert _{1}\sum_{i=1}^{n}\frac{p_{i}^{2}}{%
1-2( 1-\mathrm{e}^{-p_{i}}) }:=UB_{\mathit{Po}},
\end{eqnarray*}
where $\lambda =\sum_{i=1}^{n}p_{i}$ and $\Vert \Delta
^{2}f_{\mathit{Po}(\lambda )}\Vert _{1}$ is given in Proposition \ref{np}.
\end{corollary}

\begin{remark}\label{r2} If we assume that $\sum_{i=1}^{n}p_{i}^{2}\rightarrow 0$
as $n\rightarrow \infty $ (implying that $\max_{i}p_{i}\rightarrow 0$), the
first term $(\sum_{i=1}^{n}p_{i}^{2})^{2}$ in the upper bound $\mathit{UB}_{\mathit{Po}}$
(Corollary \ref{iP}) or in $\mathit{UB}_{\mathit{CP}}$ (Theorem \ref{iCP}) tends to $0$ at a
faster rate than the second term and, therefore, the order of $\mathit{UB}_{\mathit{Po}}$ and $%
\mathit{UB}_{\mathit{CP}}$ is the same as the order of their second term.  That is, for $\mathit{UB}_{\mathit{CP}}$,
we have
\begin{eqnarray*}
\mathit{UB}_{\mathit{CP}}\sim
\cases{\displaystyle
{1/4}\bigl\Vert \Delta ^{2}f_{\mathit{CP}(\lambda ,F)}\bigr\Vert
_{1}\sum\limits_{i=1}^{n}E(X_{i})^{2}, & \quad \mbox{for }$\lambda  \mbox{ fixed,%
}$ \cr
\displaystyle \dfrac{1}{\lambda \mu _{2}\sqrt{2\curpi \mathrm{e}}}\sum\limits_{i=1}^{n}E(X_{i})^{2},
&\quad \mbox{when }$\lambda \rightarrow \infty $,
}
\end{eqnarray*}
where $\mu _{2}$ denotes the second  moment of the compounding
distribution $F$ (see Remark \ref{r1} above). According to Remark~\ref{r1}, the
second asymptotic result for $\mathit{UB}_{\mathit{CP}}$ above (when \mbox{$\lambda \rightarrow \infty $})
is valid when $\mathit{CP}(\lambda ,F)$ \textit{is close to a normal distribution}.
Therefore, we can say that, for independent $X_{1},\ldots,X_{n}\in \mathbb{Z}%
_{+}$ with $E(X_{i})=\mathrm{O}(p_{i})$,
\begin{eqnarray*}
d_{\mathrm{TV}}\Biggl(\mathcal{L}\Biggl(\sum_{i=1}^{n}X_{i}\Biggr),\mathit{CP}(\lambda ,F)\Biggr)=\mathrm{O}\Biggl( \frac{1}{%
\lambda }\sum_{i=1}^{n}p_{i}^{2}\Biggr),
\end{eqnarray*}
whenever $F$ is such that $\mathit{CP}(\lambda ,F)\approx N(\mu ,\sigma ^{2})$ or,
more generally, whenever $\Vert \Delta ^{2}f_{\mathit{CP}(\lambda
,F)}\Vert _{1}=\mathrm{O}(\lambda ^{-1})$. Our approach requires the
restriction $\sum_{i=1}^{n}p_{i}^{2}\rightarrow 0$ (not only $%
\max_{i}p_{i}\rightarrow 0), $but we have reasons to believe (see Remark~\ref{r3}
) that this restriction is superfluous and can be weakened. This offers a
clue to a question raised by Le Cam (\citeyear{1960Le}) (see also Barbour and Utev (\citeyear{1999Barbour})
and Roos (\citeyear{2003Roos})) about the form of the compounding distribution $F$ that
would permit us to achieve a compound Poisson approximation error order
similar to that obtained for Poisson approximation, that is, $\frac{1}{\lambda }%
\sum_{i=1}^{n}p_{i}^{2}$.

We also point out that the upper bound $\mathit{UB}_{\mathit{Po}}$ of Corollary \ref{iP} for
the Poisson approximation is similar to the one derived by Deheuvels and
Pfeifer (\citeyear{1986Deheuvels}), (see also Deheuvels, Pfeifer and Puri (\citeyear{1989Deheuvels})) who employed
an entirely different method. The factor $\Vert \Delta
^{2}f_{\mathit{Po}(\lambda )}\Vert _{1}/4$ appears in the bounds of
these articles (in an equivalent form, not recognized as being the
$L_{1}$-norm of $\Delta ^{2}f_{\mathit{Po}(\lambda )}/4$) and was proven to
be optimal (that is, $d_{\mathrm{TV}}\sim $ $\mathit{UB}_{\mathit{Po}}$; see Deheuvels and Pfeifer
(\citeyear{1986Deheuvels})) under the usual asymptotic assumptions. The same argument is
possibly true for the more general smoothness factor $\Vert \Delta
^{2}f_{\mathit{CP}(\lambda ,F)}\Vert _{1}/4$.
\end{remark}

\begin{remark}\label{r3} In the proof of Theorem \ref{iCP}, the quantity $%
2\sum_{m=1}^{n}a_{m}b_{m}$ (see relation (\ref{eq23})) was bounded rather
crudely in order to obtain a closed form upper bound. This resulted in a
simple-in-form first term, namely $(\sum_{i=1}^{n}p_{i}^{2})^{2}$, in $%
\mathit{UB}_{\mathit{CP}}$. If $\sum_{i=1}^{n}p_{i}^{2}\rightarrow 0$, then this term does not
have a significant effect on $\mathit{UB}_{\mathit{CP}}$, but if $\sum_{i=1}^{n}p_{i}^{2}$ is
not close to $0$, then it may result in a very crude upper bound.
\end{remark}

Nevertheless, concerning the Poisson case, if we possessed a simple-in-form
upper bound for $\Vert \Delta ^{2}f_{\mathit{Po}(\lambda )}\Vert _{1}$, we
could obtain a better (smaller) bound for the quantity $%
2\sum_{m=1}^{n}a_{m}b_{m}$. To get an idea of how this can be done, we shall
treat  the simplest case where $X_{1},X_{2},\ldots,X_{n}$ are i.i.d. ($p_{i}=p$)
Bernoulli random variables. Recall that, in general, $\frac{1}{4}%
\Vert \Delta ^{2}f_{\mathit{Po}(\lambda )}\Vert _{1}\leq ( 1\wedge
\frac{1}{3\lambda }) $ and, therefore,
\begin{eqnarray*}
a_{m} &\leq &{1/2}\Vert \Delta ^{2}f_{\Sigma
_{i=m+1}^{n}N_{i}}\Vert _{1}\mathbf{\zeta }_{2}(X_{m},N_{m})\leq
\biggl( 1\wedge \frac{1}{3(n-m)p}\biggr) p^{2},
\\
b_{m} &\leq &\frac{1-\mathrm{e}^{-(m-1)p}}{(m-1)p}\sum_{i=1}^{m-1}p_{i}^{2}\leq p.
\end{eqnarray*}

 Assuming that $\lambda \geq 1/3+p$ and taking into account that $%
\sum_{i=n_{1}}^{n_{2}}\frac{1}{i}<\log ( \frac{n_{2}}{n_{1}-1}%
) $, the sum $2\sum_{m=1}^{n}a_{m}b_{m}$ is bounded above by
\begin{eqnarray}\label{eq49}
2\sum_{m=1}^{n}\biggl( 1\wedge \frac{1}{3(n-m)p}\biggr) p^{3} &\leq
&\sum_{m=1}^{\lfloor n-{1/(3p)}\rfloor }\frac{2p^{2}}{3(n-m)}%
+\sum_{m=\lfloor n-{1/(3p)}\rfloor +1}^{n}2p^{3}\nonumber
 \\[-8pt]\\[-8pt]
&\leq & 2\frac{p^{2}}{3}\biggl( \log \frac{3np}{1-3p}+1\biggr) +2p^{3},\nonumber
\end{eqnarray}
under the assumption $p<1/3$. The latter reveals that, when $p_{i}=p$ and $%
\lambda >1/3+p$, the term $(\sum_{i=1}^{n}p_{i}^{2})^{2}=\lambda ^{2}p^{2}$
in Corollary \ref{iP} can be substantially reduced to (\ref{eq49}), implying
that $\mathit{UB}_{\mathit{Po}}\approx \frac{2}{3}p^{2}(\log 3\lambda +1)+\frac{1}{\sqrt{2\curpi \mathrm{e}%
}}p$. The above bound could also be reduced (requiring more complicated
algebraic manipulations) in the case of non-i.i.d. Bernoulli random
variables. For a more general case though, for example, in a compound Poisson
approximation, we must first find a suitable general upper bound for $%
\Vert \Delta ^{2}f_{\mathit{CP}(\lambda ,F)}\Vert _{1}$ which, at the
moment, does not seem an easy task and is left for future work.

%s4 ###
\section{Compound Poisson approximation for sums of $k$-dependent random
variables}\label{sec4}

In this section, a more general setup is assumed. We are now interested in
approximating the distribution of the sum $X_{1}+\cdots+X_{n}$ when the $k$%
-dependent $X_{i}$'s are rarely non-zero. Naturally, we expect that this
distribution converges weakly to an appropriate compound Poisson
distribution.

Following the same methodological steps as in the proof of Theorem \ref{iCP},
we offer a bound that includes a smoothness factor analogous to a Stein
factor. The appearance of such a factor is perhaps  the first  (for sums of
\textit{dependent} random variables) outside the Stein--Chen method. As
was mentioned in the \hyperref[sec1]{Introduction}, the smoothness factor we derive is
simpler, seems more natural and is better than the corresponding Stein
factors. On the other hand, inevitably, an undesired term analogous to $%
(\sum p_{i}^{2})^{2}$ of Theorem \ref{iCP} again appears in the upper bounds.

For convenience, we shall focus our approach on a sequence of independent
random variables $Z_{1},Z_{2},\ldots$ defined over a probability space $(\Omega
,\mathcal{A},P)$ and consider $k$-dependent random variables of the form $%
h_{i}(Z_{i},\ldots,Z_{i+k-1})$. This approach is not restrictive since, in
almost all applications, local dependency arises in this setup (for example, runs
or scan statistics, patterns, reliability theory, graph theory problems,
moving sums, etcetera). Specifically, let $Z_{1},Z_{2},\ldots$ be independent random
variables and  also let
%
%e13 ###
\begin{equation}\label{eq26}
X_{i}=h_{i}(Z_{i},\ldots,Z_{i+k-1}),\qquad i=1,2,\ldots,
\end{equation}
be a sequence of non-negative, integer-valued random variables,
generated by some measurable functions $h_{i}\dvtx \mathbb{R}^{k}\rightarrow
\mathbb{Z}$. The above definition implies that $X_{i}$ is independent of $%
X_{1},\ldots,X_{i-k}$ and $X_{i+k},\ldots.$ Therefore, $X_{1},X_{2},\ldots$ are
``$k$-dependent'' random variables
(independent random variables can be considered as $1$-dependent).
Naturally, the bound offered tends to $0$, provided that $P(X_{i}$ $\neq
0)=p_{i}\approx 0$. Hence, the condition $\max_{i}\sum_{j=i-3k+3}^{i}p_{j}<%
\log 2\approx 0.693$ does not affect the generality of the result. We assume
that $X_{i}=0$ for all $i<1$.

\begin{theorem}
\label{kCP}Let $X_{1},X_{2},\ldots,X_{n}\in \mathbb{Z}_{+}$ be $k$-dependent
random variables (defined as in (\ref{eq26})) with finite second moments.
Let $N_{1},\ldots ,N_{n}$ be independent random variables (also independent of $%
Z_{i}$) with $N_{i}$ following the $\mathit{CP}(p_{i},G_{i})$ distribution, where $%
G_{i}(x)=P(X_{i}\leq x|X_{i}\neq 0)$, $x\in \mathbb{R}$ and $%
 p_{i}=P(X_{i}\neq 0)$. Then, for $m:=\max_{i}\sum_{j=i-3k+3}^{i}p_{j}<\log
2 $,
\begin{eqnarray*}
&& d_{\mathrm{TV}}\Biggl( \mathcal{L}\sum_{i=1}^{n}X_{i},\mathit{CP}(\lambda
_{n},F_{n})\Biggr)
\\
&&\quad \leq C_{n}+\frac{\Vert \Delta ^{2}f_{\mathit{CP}(\lambda
_{n},F_{n})}\Vert _{1}}{2( 1-2(1-\mathrm{e}^{-m})) }\sum_{i=1}^{n}%
\mathbf{\zeta }_{2}\Biggl( \mathcal{L}\sum\limits_{j=i-2k+2}^{i}X_{j},%
\mathcal{L}\Biggl( \sum\limits_{j=i-2k+2}^{i-1}X_{j}+N_{i}\Biggr) \Biggr)
\\
&&\quad :=\mathit{UB}_{\mathit{CP}}^{\prime },
\end{eqnarray*}
where
\begin{eqnarray*}
C_{n}&:= &2\sum_{i=1}^{n}\Biggl( d_{\mathrm{TV}}\Biggl( \mathcal{L}%
\sum_{j=1}^{i-3k+2}X_{j},\mathcal{L}\sum_{j=1}^{i-3k+2}N_{j}\Biggr)
+\sum_{j=i-3k+3}^{i-2k+1}p_{j}\Biggr)
\\
&&{}\qquad \times \Biggl( 2P\bigl( (X_{i-k+1},\ldots,X_{i-1})\neq \mathbf{0},X_{i}\neq
0\bigr) +2p_{i}\sum_{j=i-k+1}^{i-1}p_{j}+p_{i}^{2}\Biggr)
\end{eqnarray*}
and $\lambda_{n}=\sum_{i=1}^{n}p_{i}$, $F_{n}=\sum_{i=1}^{n}\frac{p_{i}}{%
\lambda_{n}}G_{i}$.
\end{theorem}

\begin{pf}
In order to simplify notation, set $X_{a,b}:=\sum_{i=a}^{b}X_{i}$, $\mathbf{%
X}_{a,b}:=(X_{a},X_{a+1},\ldots,X_{b})$, $N_{a,b}:=\sum_{i=a}^{b}N_{i}$ and $%
\mathbf{Z}_{a,b}:=(Z_{a},Z_{a+1},\ldots,Z_{b})$. Also, let $U_{i},U_{i}^{\ast }$,
$i=1,2,\ldots,n$, be independent random variables, also independent of $%
Z_{i},N_{i}$ following the uniform distribution on $(0,1)$.

Fix $i\in \{1,2,\ldots,n\}$. In order to avoid a special treatment for small
values of $i$ due to edge effects and to preserve a unified analysis for all $%
i $ that takes into account edge effects, we simply assume that $%
X_{j}=N_{j}=Z_{j}=0$ for $j\leq 0$. According to Lemma \ref{lemma6}(b)
(with $f$ being the p.d.f. of $N_{1,i-3k+2}$), there exists a random
variable $N_{1,i-3k+2}^{\ast }=g_{1}(U_{i-3k+2},X_{1,i-3k+2})$ such that $%
\mathcal{L}N_{1,i-3k+2}^{\ast }=\mathcal{L}N_{1,i-3k+2}$ and $%
(X_{1,i-3k+2},N_{1,i-3k+2}^{\ast})$ are maximally coupled, that is,
\begin{eqnarray*}
d_{\mathrm{TV}}(\mathcal{L}X_{1,i-3k+2},\mathcal{L}N_{1,i-3k+2}^{\ast})
=P(X_{1,i-3k+2}\neq N_{1,i-3k+2}^{\ast }).
\end{eqnarray*}
Moreover, according to Lemma \ref{lemma6}(a)
(with $f$ now being the p.d.f. of $N_{i}$), there exists a random variable
\begin{eqnarray*}
N_{i}^{\ast }=g_{2}(U_{i}^{\ast },X_{i},\mathbf{X}_{i-k+1,i-1},\mathbf{Z}%
_{i-k+1,i-1})
\end{eqnarray*}
such that $\mathcal{L}N_{i}^{\ast }=\mathcal{L}N_{i}$, $N_{i}^{\ast }$ is
independent of the vector $(\mathbf{X}_{i-k+1,i-1},\mathbf{Z}_{i-k+1,i-1})$
and
%
%e14 ###
\begin{equation}\label{eq6}
d_{\mathrm{TV}}(\mathcal{L}(X_{i},\mathbf{X}_{i-k+1,i-1},\mathbf{Z}_{i-k+1,i-1}),%
\mathcal{L}(N_{i}^{\ast },\mathbf{X}_{i-k+1,i-1},\mathbf{Z}%
_{i-k+1,i-1}))=P( X_{i}\neq N_{i}^{\ast }).
\end{equation}
It is easy to check that, as defined, $N_{i}^{\ast }$ is also independent of
$\mathbf{X}_{1,i-1}$. Indeed, if we set $\mathbf{Y:}=(\mathbf{Z}_{i-k+1,i-1},%
\mathbf{X}_{i-k+1,i-1})$, for all $x,\mathbf{x}$, we have that%
\begin{eqnarray*}
P( N_{i}^{\ast }=x,\mathbf{X}_{1,i-1}=\mathbf{x}) =\sum_{\mathbf{y%
}}P(N_{i}^{\ast }=x,\mathbf{Y}=\mathbf{y},\mathbf{X}_{1,i-1}=\mathbf{x}).
\end{eqnarray*}
We may write $\mathbf{X}_{1,i-1}=\mathbf{g}(\mathbf{Z}_{1,i-k},\mathbf{Y})$
for some appropriate function $\mathbf{g}$ taking values in $\mathbb{Z}%
^{i-1} $. Hence, the above sum is equal to
\begin{eqnarray*}
&& \sum_{\mathbf{y}}P\bigl( g_{2}(U_{i}^{\ast },X_{i},\mathbf{Y}%
)=x,\mathbf{Y=y,g}(\mathbf{Z}_{1,i-k},\mathbf{y})=\mathbf{x}\bigr)
 \\
&&\quad =\sum_{\mathbf{y}}P\bigl(g_{2}(U_{i}^{\ast },X_{i},\mathbf{Y})=x,\mathbf{Y=y}\bigr)P\bigl(%
\mathbf{g}(\mathbf{Z}_{1,i-k},\mathbf{y})=\mathbf{x}\bigr)
\\
&&\quad =\sum_{\mathbf{y}}P(N_{i}^{\ast }=x)P(\mathbf{Y}=\mathbf{y})P\bigl(\mathbf{g}(%
\mathbf{Z}_{1,i-k},\mathbf{y})=\mathbf{x}\bigr)
 \\
&&\quad =\sum_{\mathbf{y}}P(N_{i}^{\ast }=x)P\bigl(\mathbf{Y}=\mathbf{y},\mathbf{g}(%
\mathbf{Z}_{1,i-k},\mathbf{y})=\mathbf{x}\bigr)
\\
&&\quad =P(N_{i}^{\ast }=x)\sum_{\mathbf{y}}P( \mathbf{Y}=\mathbf{y,X}%
_{1,i-1}=\mathbf{x}) =P(N_{i}^{\ast }=x)P(\mathbf{X}_{1,i-1}=\mathbf{y}%
),
\end{eqnarray*}
which is valid for all $x,\mathbf{y}$ and thus $N_{i}^{\ast }$ is
independent of $\mathbf{X}_{1,i-1}$.

Now, applying the inequality (see Lemma \ref{lemma3})
\begin{eqnarray*}
d_{\mathrm{TV}}\bigl(\mathcal{L}(Z+X),\mathcal{L}(Z+Y)\bigr)\leq 2P(Z\neq W,X\neq Y)+d_{\mathrm{TV}}\bigl(%
\mathcal{L}(W+X),\mathcal{L}(W+Y)\bigr)
\end{eqnarray*}
with $Z=X_{1,i-1},X=X_{i}+N_{i+1,n},Y=N_{i}^{\ast
}+N_{i+1,n},W=N_{1,i-3k+2}^{\ast }+X_{i-2k+2,i-1}$, we obtain
\begin{eqnarray} \label{eq7}
&& d_{\mathrm{TV}}\bigl(\mathcal{L}(X_{1,i-1}+X_{i}+N_{i+1,n}),\mathcal{L}%
(X_{1,i-1}+N_{i}^{\ast }+N_{i+1,n})\bigr)\nonumber
 \\
&&\quad \leq 2P(X_{1,i-1}\neq N_{1,i-3k+2}^{\ast
}+X_{i-2k+2,i-1},X_{i}+N_{i+1,n}\neq N_{i}^{\ast }+N_{i+1,n})\nonumber
\\[-8pt]\\[-8pt]
&&\qquad{}+d_{\mathrm{TV}}\bigl(\mathcal{L}(N_{1,i-3k+2}^{\ast }+X_{i-2k+2,i}+N_{i+1,n}),\nonumber
\\
&&\qquad \qquad\quad\,\mathcal{L%
}(N_{1,i-3k+2}^{\ast }+X_{i-2k+2,i-1}+N_{i}^{\ast }+N_{i+1,n})\bigr).\nonumber
\end{eqnarray}
Note that $\mathcal{L}N_{i}^{\ast }=\mathcal{L}N_{i}$ and $N_{i}^{\ast }$ is
independent of $\mathbf{X}_{1,i-1}$, also that $\mathcal{L}N_{1,i-3k+2}^{\ast
}=\mathcal{L}N_{1,i-3k+2}$ and $N_{1,i-3k+2}^{\ast
}=g_{1}(U_{i-3k+2},X_{1,i-3k+2})$ is independent of $\mathbf{X}_{i-2k+2,i}$
and $N_{i}^{\ast }$. Therefore, we have that
\begin{eqnarray*}
\mathcal{L}( X_{1,i-1}+N_{i}^{\ast }+N_{i+1,n}) &=&\mathcal{L}
( X_{1,i-1}+N_{i}+N_{i+1,n}) ,
\\
\mathcal{L}(N_{1,i-3k+2}^{\ast }+X_{i-2k+2,i}+N_{i+1,n})&=&\mathcal{L}
(N_{1,i-3k+2}+X_{i-2k+2,i}+N_{i+1,n}),
\\
\mathcal{L}(N_{1,i-3k+2}^{\ast }+X_{i-2k+2,i-1}+N_{i}^{\ast }+N_{i+1,n})&=&
\mathcal{L}(N_{1,i-3k+2}+X_{i-2k+2,i-1}+N_{i}+N_{i+1,n}).
\end{eqnarray*}
Using the above relations, inequality (\ref{eq7}) is equivalent to
%
%e15 ###
\begin{equation} \label{eq27}
d_{\mathrm{TV}}\bigl(\mathcal{L}(X_{1,i-1}+X_{i}+N_{i+1,n}),\mathcal{L}%
(X_{1,i-1}+N_{i}+N_{i+1,n})\bigr)\leq 2a_{i}+b_{i},
\end{equation}
where
\begin{eqnarray*}
a_{i} &=&P(X_{1,i-2k+1}\neq N_{1,i-3k+2}^{\ast },X_{i}\neq N_{i}^{\ast }),
\\
b_{i} &=&d_{\mathrm{TV}}\bigl(\mathcal{L}(N_{1,i-3k+2}+X_{i-2k+2,i-1}+X_{i}+N_{i+1,n}),
\\
&&\qquad \mathcal{L}(N_{1,i-3k+2}+X_{i-2k+2,i-1}+N_{i}+N_{i+1,n})\bigr).
\end{eqnarray*}
The random variables $X_{i},N_{i}^{\ast }$ are independent of $%
X_{1},\ldots,X_{i-2k+1},N_{1,i-3k+2}^{\ast }$ and, hence, it is easy to see that
%
%e16 ###
\begin{eqnarray}\label{eq28}
a_{i} &=&P( X_{1,i-3k+2}+X_{i-3k+3,i-2k+1}\neq N_{1,i-3k+2}^{\ast
})P(X_{i}\neq N_{i}^{\ast })\nonumber
 \\
&\leq &\bigl( P(X_{1,i-3k+2}\neq N_{1,i-3k+2}^{\ast })+P\mathbb{(}%
X_{i-3k+3,i-2k+1}\neq 0)\bigr) P(X_{i}\neq N_{i}^{\ast })
 \\
&\leq &\Biggl( d_{\mathrm{TV}}(\mathcal{L}X_{1,i-3k+2},\mathcal{L}N_{1,i-3k+2})+%
\sum_{j=i-3k+3}^{i-2k+1}p_{j}\Biggr) P(X_{i}\neq N_{i}^{\ast }).\nonumber
\end{eqnarray}
Using relation (\ref{eq6}) above along with $\mathcal{L}(\mathbf{Z}%
_{i-k+1,i-1}\mathbf{,X}_{i-k+1,i-1},N_{i}^{\ast })=\mathcal{L}(\mathbf{Z}%
_{i-k+1,i-1},\break \mathbf{X}_{i-k+1,i-1},N_{i})$, we observe that
\begin{eqnarray*}
P( X_{i}\neq N_{i}^{\ast }) =d_{\mathrm{TV}}(\mathcal{L}(\mathbf{Z}%
_{i-k+1,i-1}\mathbf{,X}_{i-k+1,i-1},X_{i}),\mathcal{L}(\mathbf{Z}_{i-k+1,i-1}%
\mathbf{,X}_{i-k+1,i-1},N_{i}))
\end{eqnarray*}
and applying Lemma \ref{lemma3} with $\mathbf{X}=X_{i}\mathbf{,Y}=N_{i}%
\mathbf{,Z}=\mathbf{(Z}_{i-k+1,i-1}\mathbf{,X}_{i-k+1,i-1}), \mathbf{W}=\break (%
\mathbf{Z}_{i-k+1,i-1}\mathbf{,0)}$, we deduce
\begin{eqnarray} \label{eq29}
&&P( X_{i}\neq N_{i}^{\ast })\nonumber
 \\
&&\quad \leq 2P(X_{i}\neq N_{i},\mathbf{X}_{i-k+1,i-1}\neq \mathbf{0})+d_{\mathrm{TV}}(%
\mathcal{L}(\mathbf{Z}_{i-k+1,i-1}\mathbf{,0},X_{i}),\mathcal{L}(\mathbf{Z}%
_{i-k+1,i-1}\mathbf{,0},N_{i}))\nonumber
 \\
&&\quad \leq 2P(X_{i}\neq 0,\mathbf{X}_{i-k+1,i-1}\neq \mathbf{0})\nonumber
\\[-8pt]\\[-8pt]
&&\qquad{}+2P(N_{i}\neq
0)P(\mathbf{X}_{i-k+1,i-1}\neq \mathbf{0})+d_{\mathrm{TV}}(\mathcal{L}X_{i},\mathcal{L%
}N_{i})\nonumber
 \\
&&\quad \leq 2P(X_{i}\neq 0,\mathbf{X}_{i-k+1,i-1}\neq \mathbf{0}%
)+2p_{i}\sum_{j=i-k+1}^{i-1}p_{j}+p_{i}^{2}.\nonumber
\end{eqnarray}
Next, we consider a random variable $N^{\perp }$ with $\mathcal{L}%
N_{i}^{\perp }=\mathcal{L}N_{i}$,  independent of all other random
variables involved in our analysis. Applying the inequality (\ref{eq17})
with $W=N_{i-3k+3,i}$ and assuming that $P(N_{i-3k+3,i}\neq 0)<1/2$ (which
is valid since we have assumed that $m<\log 2)$, we get
\begin{eqnarray*}
b_{i} &=& d_{\mathrm{TV}}\bigl(\mathcal{L}(N_{1,i-3k+2}+X_{i-2k+2,i-1}+X_{i}+N_{i+1,n}),%
\\
&&\qquad\mathcal{L}(N_{1,i-3k+2}+X_{i-2k+2,i-1}+N_{i}^{\perp }+N_{i+1,n})\bigr)
\\
&\leq &\bigl( 1-2P(N_{i-3k+3,i}\neq 0)\bigr) ^{-1} d_{\mathrm{TV}}\bigl(
\mathcal{L}(N_{1,n}+X_{i-2k+2,i}),\mathcal{L}(N_{1,n}+X_{i-2k+2,i-1}+N_{i}^{%
\perp })\bigr).
\end{eqnarray*}
Finally, using Lemma \ref{Btv}, we derive
%
%e17 ###
\begin{equation}\label{eq30}
b_{i}\leq \frac{{1/2}\Vert \Delta ^{2}f_{\mathit{CP}(\lambda
,F)}\Vert _{1}}{1-2(1-\mathrm{e}^{-\sum_{j=i-3k+3}^{i}p_{j}})}\mathbf{\zeta }%
_{2}\bigl(\mathcal{L}X_{i-2k+2,i},\mathcal{L}(X_{i-2k+2,i-1}+N_{i})\bigr).
\end{equation}
Combining (\ref{eq27})--(\ref{eq29}) with (\ref{eq30}), we
obtain, for all $i=1,2,\ldots,n$, the  inequality
\begin{eqnarray*}
&& d_{\mathrm{TV}}\bigl(\mathcal{L}(X_{1,i-1}+X_{i}+N_{i+1,n}),\mathcal{L}%
(X_{1,i-1}+N_{i}+N_{i+1,n})\bigr)
\\
&&\quad\leq 2\Biggl( d_{\mathrm{TV}}(\mathcal{L}X_{1,i-3k+2},\mathcal{L}N_{1,i-3k+2})+%
\sum_{j=i-3k+3}^{i-2k+1}p_{j}\Biggr)
\\
&&\qquad {}\times \Biggl( 2P(X_{i}\neq 0,\mathbf{X}_{i-k+1,i-1}\neq \mathbf{0}%
)+2p_{i}\sum_{j=i-k+1}^{i-1}p_{j}+p_{i}^{2}\Biggr)
\\
&&\qquad{} +\frac{{1/2}\Vert \Delta ^{2}f_{\mathit{CP}(\lambda ,F)}\Vert _{1}%
}{1-2(1-\mathrm{e}^{-m})}\mathbf{\zeta }_{2}\bigl(\mathcal{L}X_{i-2k+2,i},\mathcal{L}%
(X_{i-2k+2,i-1}+N_{i})\bigr)
\end{eqnarray*}
and the final result follows immediately by virtue of the Lindeberg
decomposition (triangle inequality)
\begin{eqnarray*}
\quad d_{\mathrm{\mathrm{TV}}}( \mathcal{L}X_{1,n},\mathcal{L}N_{1,n}) \leq
\sum_{i=1}^{n}d_{\mathrm{TV}}\bigl(\mathcal{L}(X_{1,i-1}+X_{i}+N_{i+1,n}),\mathcal{L}%
(X_{1,i-1}+N_{i}+N_{i+1,n})\bigr).\quad\qed
\end{eqnarray*}
\noqed\end{pf}

\begin{remark}\label{r4} The upper bound $\mathit{UB}_{\mathit{CP}}^{\prime }$ in Theorem \ref{kCP}
is composed of two terms, the first of which is the quantity $C_{n}$,
which is analogous to the term $(\Sigma p_{i}^{2})^{2}$ appearing in Theorem %
\ref{iCP} that concerns the independent summands case. As it was for $%
(\Sigma p_{i}^{2})^{2}$, the term $C_{n}$ tends to $0$ faster than the
second term of $\mathit{UB}_{\mathit{CP}}^{\prime }$, under certain asymptotic conditions.
Therefore, under these conditions, the order of $\mathit{UB}_{\mathit{CP}}^{\prime }$
coincides with the order of the second term.
\end{remark}

\begin{remark}\label{r5} If $X_{1},X_{2},\ldots ,X_{n}$ are $k$-dependent Bernoulli
random variables then, similarly to Corollary \ref{iP}, Theorem \ref{kCP}
implies a Poisson approximation result. Specifically, Theorem \ref{kCP} can
now be written with $\mathit{Po}(p_{i})$ in place of $\mathit{CP}(p_{i},G_{i})$ and $%
\mathit{Po}(\lambda _{n})$ in place of $\mathit{CP}(\lambda _{n},F_{n})$. Consequently, the
norm $\Vert \Delta ^{2}f_{\mathit{Po}(\lambda _{n})}\Vert _{1}$ will
appear instead of $\Vert \Delta ^{2}f_{\mathit{CP}(\lambda
_{n},F_{n})}\Vert _{1}$.
\end{remark}\vspace{-2pt}

The upper bound $\mathit{UB}_{\mathit{CP}}^{\prime }$ in Theorem \ref{kCP} may seem difficult
to apply in its present form. For this reason, we present the following
corollary  which provides two slightly worse, but more easily computable, upper
bounds. The bound (a) is valid without any assumption on the form of
dependence among  the $X_{i}$'s. The bound (b) is smaller than (a), but is valid
only when the $X_{i}$'s exhibit a certain weak form of positive/negative
dependence. We recall that two random variables $Y_{1},Y_{2}$ are called
\textit{positively quadrant dependent }(PQD) if
%
%e18 ###
\begin{equation}\label{eq32}
P(Y_{1}\geq x_{1},Y_{2}\geq x_{2})\geq P(Y_{1}\geq x_{1})P(Y_{2}\geq
x_{2}) \qquad \mbox{for all }  x_{1},x_{2}
\end{equation}
and \textit{negatively quadrant dependent} (\textit{NQD}) if (\ref{eq32})
holds, but with the inequality sign reversed. Manifestly, if $%
X_{1},X_{2},\ldots ,X_{n}$ are associated (resp., negatively associated), then
the random variables $X_{j}+X_{2}+\cdots +X_{i-1}$ and $X_{i}$ are PQD
(resp., NQD) for every $1\leq j<i\leq n$. Therefore, part (b) of the next
corollary remains valid under the stronger condition of association or
negative association of $X_{i}$'s.\vspace{-2pt}

\begin{corollary}
\label{cor3}\textup{(a)} Let $X_{1},X_{2},\ldots,X_{n}\in \mathbb{Z}_{+}$ be $k$%
-dependent random variables (defined as in (\ref{eq26})) with finite second
moments. Then, for $m:=\max_{i}\sum_{j=i-3k+3}^{i}p_{j}<\log 2,$ $%
p_{i}=P(X_{i}\neq 0)$,
\begin{eqnarray*}
&& d_{\mathrm{TV}}\Biggl( \mathcal{L}\sum_{i=1}^{n}X_{i},\mathit{CP}(\lambda
_{n},F_{n})\Biggr)
\\[-2pt]
&&\quad \leq C_{n}+\frac{\Vert \Delta ^{2}f_{\mathit{CP}(\lambda
_{n},F_{n})}\Vert _{1}}{2( 1-2(1-\mathrm{e}^{-m})) }%
\sum_{i=1}^{n}\Biggl( \sum_{j=i-k+1}^{i-1}\bigl(
E(X_{i}X_{j})+E(X_{i})E(X_{j})\bigr) +\frac{1}{2}E(X_{i})^{2}\Biggr)
\\[-3pt]
&&\quad :=\mathit{UB}_{\mathit{CP}}^{\prime \prime },
\end{eqnarray*}\vspace{-3pt}
where\vspace{-3pt}
\begin{eqnarray*}
C_{n} &:=& 2\sum_{i=1}^{n}\Biggl( 2\sum_{j=1}^{i-3k+2}\Biggl(
\sum_{t=j-k+1}^{j-1}\bigl( P( X_{t}X_{j}\neq 0)
+p_{t}p_{j}\bigr) +\frac{1}{2}p_{j}^{2}\Biggr)
+\sum_{j=i-3k+3}^{i-2k+1}p_{j}\Biggr)
\\[-2pt]
&&\times \Biggl( 2\sum_{j=i-k+1}^{i-1}\bigl( P(X_{j}\neq 0,X_{i}\neq
0)+p_{i}p_{j}\bigr) +p_{i}^{2}\Biggr)
\end{eqnarray*}
and $\lambda _{n}=\sum_{i=1}^{n}p_{i},$ $F_{n}=\sum_{i=1}^{n}\frac{p_{i}}{%
\lambda _{n}}G_{i},G_{i}(x)=P(X_{i}\leq x|X_{i}\neq 0),x\in \mathbb{R}$ ($%
X_{i}=0$ for all $i<1$).

\textup{(b)} If, in addition, the random variables $X_{j}+\cdots+X_{i-1}$ and $X_{i}$
are PQD or NQD for every $1\leq j<i\leq n$, then the bound $\mathit{UB}_{\mathit{CP}}^{\prime
\prime }$ in \textup{(a)} is valid with $\vert \operatorname{Cov}(X_{i},X_{j})\vert $ in
place of $E(X_{i}X_{j})+E(X_{i})E(X_{j})$ and
\begin{eqnarray*}
\sum_{t=j-k+1}^{j-1}\vert \operatorname{Cov}(X_{j},X_{t})\vert +\frac{1}{2}%
E(X_{j})^{2} \quad  \mbox{in place of}\quad \sum_{t=j-k+1}^{j-1}\bigl( P(
X_{t}X_{j}\neq 0) +p_{t}p_{j}\bigr) +\frac{1}{2}p_{j}^{2}.
\end{eqnarray*}
\end{corollary}

\begin{pf}
\textup{(a)}  This follows readily from Theorem \ref{kCP} by applying Corollary \ref%
{cor4} above, Corollary 7 in Boutsikas (\citeyear{2006Boutsikas}) and the fact that
\begin{eqnarray*}
\mathbf{\zeta }_{2}\Biggl( \mathcal{L}\Biggl(
\sum_{j=l}^{i-1}X_{j}+X_{i}^{\perp }\Biggr) ,\mathcal{L}\Biggl(
\sum_{j=l}^{i-1}X_{j}+N_{i}\Biggr) \Biggr) \leq \mathbf{\zeta }_{2}(
\mathcal{L}X_{i}^{\perp },\mathcal{L}N_{i}) =\frac{1}{2}E(X_{i})^{2}
\end{eqnarray*}
($N_{i}\sim \mathit{CP}(p_{i},G_{i})), $which is a consequence of the regularity
property of \textbf{$\zeta $}$_{2}$ combined with equality (\ref{eq36}).

\textup{(b)} This is again a direct consequence of Theorem \ref{kCP}. Set $%
W=\sum\nolimits_{j=i-2k+2}^{i-1}X_{j}$ and let $X_{i}^{\perp }$ be a random
variable independent of all other random variables involved in our analysis
with $\mathcal{L}X_{i}=\mathcal{L}X_{i}^{\perp }$. Assume that $%
X_{j}+\cdots+X_{i-1}$ and $X_{i}$ are PQD for all $j<i$. Thus, $W$ and $X_{i}$ are
PQD and hence $W+X_{i}$ is larger than $W+X_{i}^{\perp }$ with respect to
the convex order (see Section 3.3 in Boutsikas and Vaggelatou (\citeyear{2002Boutsikas})).
Therefore,\vspace{-3pt}
\begin{eqnarray*}
\mathbf{\zeta }_{2}\bigl( \mathcal{L}( W+X_{i}) ,\mathcal{L}%
( W+X_{i}^{\perp }) \bigr) =\dfrac{1}{2}\bigl(
\operatorname{Var}(W+X_{i})-\operatorname{Var}(W+X_{i}^{\perp })\bigr)
=\sum_{j=i-k+1}^{i-1}\operatorname{Cov}(X_{i},X_{j}).
\end{eqnarray*}\vspace{-3pt}
Since $W$ is independent of $X_{i}^{\perp },N_{i}$, the regularity property
of \textbf{$\zeta $}$_{2}$ and equality (\ref{eq36}) guarantee that
\begin{eqnarray*}
\mathbf{\zeta }_{2}\bigl( \mathcal{L}( W+X_{i}^{\perp }) ,%
\mathcal{L}( W+N_{i}) \bigr) \leq \mathbf{\zeta }_{2}(
\mathcal{L}X_{i}^{\perp },\mathcal{L}N_{i}) =\tfrac{1}{2}E(X_{i})^{2}.
\end{eqnarray*}
Hence, using the triangle inequality and the above two equalities, we
deduce that
\begin{eqnarray*}
\mathbf{\zeta }_{2}\bigl( \mathcal{L}( W+X_{i}) ,\mathcal{L}%
( W+N_{i}) \bigr) \leq \sum_{j=i-k+1}^{i-1}\operatorname{Cov}(X_{i},X_{j})+%
\frac{1}{2}E(X_{i})^{2}.
\end{eqnarray*}
Furthermore, from (\ref{eq18}) and Theorem 7 in Boutsikas and Vaggelatou
(\citeyear{2002Boutsikas}), we get that
\begin{eqnarray*}
&&d_{\mathrm{TV}}\Biggl( \mathcal{L}\sum_{j=1}^{i-3k+2}X_{j},\mathit{CP}(\lambda
_{i-3k+2},F_{i-3k+2})\Biggr)
\\
&&\quad \leq 2\mathbf{\zeta }_{2}\Biggl(\mathcal{L}\sum_{j=1}^{i-3k+2}X_{j},\mathit{CP}(\lambda
_{i-3k+2},F_{i-3k+2})\Biggr)
\\
&&\quad = 2\sum_{t=2}^{i-3k+2}\sum_{j=t-k+1}^{t-1}\operatorname{Cov}(X_{j},X_{t})+%
\sum_{j=1}^{i-3k+2}E(X_{j})^{2}.
\end{eqnarray*}

Similar reasoning proves the NQD random variables case (in place of
all $\operatorname{Cov}(X_{j},X_{t})$, we now get $-\operatorname{Cov}(X_{j},X_{t})>0)$.
\end{pf}

%s5 ###
\section{Illustrating applications}\label{sec5}

The purpose of this section is to illustrate the applicability and
effectiveness of the results presented in the previous sections. These
results are applicable to a wide variety of problems involving locally
dependent random variables that rarely differ from zero (for example, in risk
theory, extreme value theory, reliability theory, run and scan statistics,
graph theory and biomolecular sequence analysis). The approximation method
described in this paper, as with almost all other methods used for Poisson
approximation in the past, requires the computation of only the first- and
second-order moments of the variables involved. From this fact, it is
understood that the bounds presented  can be applied almost directly to many
of the problems where other Poisson approximation methods have been
elaborated in the past, for example, the Stein--Chen method. The main benefit of the
present method is the smoothness factor that substantially improves the
approximation error bound in many cases, while the main disadvantage is the
additional term $C_{n}$. Therefore, the conclusion here is that we usually
obtain improved bounds for moderate or small values of $\lambda$.

%s5.1 ###
\subsection{The number of overlapping runs of length $k$ in i.i.d.
trials}\label{sec5.1}

Let $\{Z_{i}\}_{i\in \mathbb{Z}}$ be a sequence of i.i.d. binary trials with
outcomes $0$ (failure) and $1$ (success), and where $P(Z_{i}=1)=p=1-q$. We are
interested in approximating the distribution of the number of (rare) \textit{%
overlapping} success runs of length $k$ within trials $1,2,\ldots,n$. This
problem has been studied in various ways by many authors in the past; see, for example,
Barbour, Holst and Janson (\citeyear{1992BarbourJan}), Balakrishnan and Koutras (\citeyear{2002Balakrishnan}) and
the references therein. We shall first derive a Poisson and then a compound
Poisson approximation.\looseness=1

(a) \textit{Poisson approximation}. If we assume that $p\rightarrow 0$ and $%
n\rightarrow \infty $, then the occurrences of success runs are rare and
asymptotically independent, and a Poisson approximation seems suitable. We
use the binary random variables
\begin{eqnarray*}
X_{i}=Z_{i}Z_{i+1}\cdots Z_{i+k-1},\qquad  i=1,2,\ldots,n-k+1.
\end{eqnarray*}%
Obviously, the random variable $\sum_{i=1}^{n-k+1}X_{i}$ counts the total
number of appearances of overlapping success runs with length $k$ which
appear within the first $n$ trials. The random variables $%
X_{1},X_{2},\ldots,X_{n-k+1}$ are $k$-dependent, can be written as in (\ref%
{eq26}) and are associated as coordinatewise non-decreasing functions of
independent random variables. Thus, they satisfy the dependence condition
required by Corollary \ref{cor3}(b). A direct application of this
corollary
for $m=(3k-2)p^{k}<\log 2$ yields

\begin{eqnarray*}
d_{\mathrm{TV}}\Biggl( \mathcal{L}\sum_{i=1}^{n-k+1}X_{i},\mathit{Po}(\lambda )\Biggr) &\leq
&C_{n-k+1}+\frac{{1/2}\Vert \Delta ^{2}f_{\mathit{Po}(\lambda
)}\Vert _{1}}{1-2(1-\mathrm{e}^{-m})}
\\
&&{}\times \Biggl( \sum_{i=2}^{n-k+1}\sum_{j=\max
\{1,i-k+1\}}^{i-1}(p^{i-j+k}-p^{2k})+\frac{n-k+1}{2}p^{2k}\Biggr)
\\
&\leq &C_{n-k+1}+\frac{\lambda p\Vert \Delta ^{2}f_{\mathit{Po}(\lambda
)}\Vert _{1}}{2q( 1-2( 1-\mathrm{e}^{-m}) ) }\biggl(
1-\biggl(k-2+\frac{1}{q}\biggr)qp^{k-1}\biggr) ,
\end{eqnarray*}
where $p_{i}=P(X_{i}=1)=p^{k}, \lambda =(n-k+1)p^{k}$ and
\begin{eqnarray*}
C_{n-k+1} &=&2\sum_{i=1}^{n-k+1}\Biggl( 2\sum_{t=2}^{i-3k+2}\sum_{j=\max
\{1,t-k+1\}}^{t-1}(p^{t-j+k}-p^{2k})+\sum_{j=1}^{i-3k+2}p^{2k}+(k-1)p^{k}%
\Biggr)
\\
&&\hspace{28pt}{}\times \Biggl( 2\sum_{j=\max \{1,i-k+1\}}^{i-1}p^{i-j+k}+(2k-1)p^{2k}\Biggr)
\\
&\leq &4\frac{\lambda ^{2}p^{2}}{q}\biggl( 1-p^{k-1}-q\biggl(k-\frac{3}{2}\biggr)p^{k-1}+%
\frac{q(k-1)p^{k-1}}{\lambda }\biggr)
\\
&&{}\times \biggl( 1-p^{k-1}+q\biggl(k-\frac{1}{2}%
\biggr)p^{k-1}\biggr).
\end{eqnarray*}
Therefore, for $m=(3k-2)p^{k}<\log 2$,
\begin{eqnarray*}
&& d_{\mathrm{TV}}\Biggl( \mathcal{L}\sum_{i=1}^{n-k+1}X_{i},\mathit{Po}(\lambda
)\Biggr)
\\
&&\quad \leq \mathit{UB}_{n,p}:=4\frac{\lambda ^{2}p^{2}}{q}\biggl( 1+\frac{qkp^{k-1}}{%
\lambda }\biggr) ( 1+qkp^{k-1}) +\frac{\lambda p\Vert \Delta
^{2}f_{\mathit{Po}(\lambda )}\Vert _{1}}{2q( 1-2( 1-\mathrm{e}^{-m})
) }.
\end{eqnarray*}
In addition, if $n\rightarrow \infty ,p\rightarrow 0$ ($k>1$ fixed), then
\begin{eqnarray*}
\mathit{UB}_{n,p}\sim
\cases{
\dfrac{\lambda \Vert \Delta ^{2}f_{\mathit{Po}(\lambda )}\Vert _{1}}{2q}p,
&\quad \mbox{when }$\lambda \mbox{ is fixed,}$ \cr
\dfrac{2}{q\sqrt{2\pi e}}p,  & \quad \mbox{when }$\lambda \rightarrow \infty $
and $p\lambda ^{2}\rightarrow 0$,%
}%
\end{eqnarray*}
where $\Vert \Delta ^{2}f_{\mathit{Po}(\lambda )}\Vert _{1}$(which is less
than $4(1\wedge \frac{1}{3\lambda })$) is given in Proposition \ref{np}. For
the same distance, a bound obtained by the Stein--Chen method (see, for example,
Barbour, Holst and Janson (\citeyear{1992BarbourJan}), page 163) is nearly equal to $2p/q$, which,
provided that $p\lambda ^{2}\approx 0$ and for moderate or large values of $%
\lambda $, is nearly four times larger than $\mathit{UB}_{n,p}$ ($\sqrt{2\curpi \mathrm{e}}%
\approx 4.1327$). For $\lambda \approx 1$, it is nearly three times
larger.

(b) \textit{Compound Poisson approximation}. The bound described in (a)
cannot help when we assume that $n\rightarrow \infty ,k\rightarrow \infty $
and $p$ is fixed. Under these conditions, the occurrences of success runs
are again rare, but they are no longer  asymptotically independent. This
happens because if a success run occurs (starts) at trial $i$ (that is, $%
Z_{i}=\cdots=Z_{i+k-1}=1$), then, with probability $p$, we shall also observe an
overlapping success run starting at position $i+1,$ and so forth. Thus, when
a success run is observed at some trial, it is likely that a number of
success runs will follow at the next trials. This ``cluster'' of adjacent
success runs is usually called a ``clump''. So, now that $n\rightarrow \infty $
and $k\rightarrow \infty $, we expect that the occurrences of clumps are
\textit{rare} and asymptotically independent, while each clump consists of
an asymptotically geometrically distributed number of overlapping success
runs. Obviously, this situation readily calls for a compound Poisson
approximation result. To achieve this, let $Y_{1},Y_{2},\ldots,Y_{n-k+1}$
represent the sizes of the clumps started at trials $1,2,\ldots,n-k+1$,
respectively. If $Y_{i}=0$, then we obviously mean that no clump has started
at position $i$. This well-known technique is called ``declumping''. More
formally, set
\begin{eqnarray*}
Y_{i}&:=&(1-Z_{i-1})\sum_{r=0}^{n-i-k+1}%
\prod_{j=i}^{i+k+r-1}Z_{j},\qquad i=2,3,\ldots,n-k+1,\quad  \mbox{and}
\\ Y_{1} &:=&
\sum_{r=0}^{n-k}\prod_{j=1}^{k+r}Z_{j}
\end{eqnarray*}
to be the size of a clump starting at position $i$ (that is, the number of adjacent
overlapping success runs until trial $n$). Clearly, $\sum_{i=1}^{n-k+1}Y_{i}$
is equal to $\sum_{i=1}^{n-k+1}X_{i}$, the total number of overlapping
success runs within trials $1,2,\ldots,n$. In this case, it is computationally
more convenient to use the \textit{stationary}, \textit{locally dependent}
random variables
\begin{eqnarray*}
Y_{i}^{\prime
}:=(1-Z_{i-1})\sum_{r=0}^{k-1}\prod_{j=i}^{i+k+r-1}Z_{j},\qquad
i=1,2,\ldots,n-k+1,
\end{eqnarray*}
which represent the \textit{truncated} sizes of clumps (their sizes cannot
be greater than $k$) starting at positions $1,2,\ldots,n-k+1$. In order to
obtain stationarity, we have also allowed the last clumps to extend further
than trial $n$. When $k$, $n$ increase so that the expected number of runs $%
(n-k+1)p^{k} $remains bounded, the processes $\mathbf{Y=(}Y_{i})\mathbf{,Y}%
^{\prime }=(Y_{i}^{\prime })$ rarely differ. This is expressed by the
following inequality (see Boutsikas (\citeyear{2006Boutsikas}), page 511):
%e19 ###
\begin{equation} \label{eq34}
d_{\mathrm{TV}}( \mathcal{L}(\mathbf{Y)},\mathcal{L}(\mathbf{Y}^{\prime
})) \leq P( \mathbf{Y}\neq \mathbf{Y}^{\prime }) \leq
(n-2k+1)qp^{2k}+2p^{k+1}.
\end{equation}
We can now use Corollary \ref{cor3}(a) to establish an upper bound for $%
d_{\mathrm{TV}}(\mathcal{L}(\sum Y_{i}^{\prime }),\mathit{CP})$. We verify that the random
variables $Y_{1}^{\prime },Y_{2}^{\prime },\ldots,Y_{n-k+1}^{\prime }\in
\mathbb{Z}_{+}$ can be written as in (\ref{eq26}) and that they are also $2k$%
-dependent. Obviously, $p_{i}=P(Y_{i}^{\prime }\neq 0)=qp^{k}$. For $%
m=(6k-2)qp^{k}<\log 2$, Corollary~\ref{cor3}(a) yields the
inequality
\begin{eqnarray*}
d_{\mathrm{TV}}\Biggl( \mathcal{L}\sum_{i=1}^{n-k+1}Y_{i}^{\prime },\mathit{CP}(\lambda
,F_{k})\Biggr) &\leq &C_{n-k+1}+\frac{\Vert \Delta ^{2}f_{\mathit{CP}(\lambda
,F_{k})}\Vert _{1}}{2( 1-2(1-\mathrm{e}^{-m})) }
\\
&&\times \sum_{i=1}^{n-k+1}\Biggl( \sum_{j=i-2k+1}^{i-1}\bigl(
E(Y_{i}^{\prime }Y_{j}^{\prime })+E(Y_{i}^{\prime })E(Y_{j}^{\prime
})\bigr) +\frac{1}{2}E(Y_{i}^{\prime })^{2}\Biggr)
\end{eqnarray*}
(we assume that $Y_{i}^{\prime }=0$ for $i<1$) with
\begin{eqnarray*}
C_{n-k+1} &\leq &2\sum_{i=1}^{n-k+1}\Biggl(
2\sum_{j=2}^{i-6k+2}\sum_{t=j-2k+1}^{j-1}\bigl(P(Y_{t}^{\prime }Y_{j}^{\prime
}\neq 0)+(qp^{k})^{2}\bigr)+\sum_{j=1}^{i-6k+2}(qp^{k})^{2}+2kqp^{k}\Biggr)
 \\
&&{}\hspace{28pt}\times \Biggl( 2\sum_{j=i-2k+1}^{i-1}P(Y_{j}^{\prime }\neq 0,Y_{i}^{\prime
}\neq 0)+4k(qp^{k})^{2}\Biggr)
\end{eqnarray*}
and $\lambda =(n-k+1)qp^{k}$, $F_{k}(x)=P(Y_{i}^{\prime }\leq
x|Y_{i}^{\prime }\neq 0)$, $x\in \mathbb{R}$. Notice  that, for $i\geq 2k$,
$P( Y_{i}^{\prime }\neq 0,Y_{j}^{\prime }\neq 0) $ is now equal to $%
q^{2}p^{2k}$ for $j=i-2k+1,\ldots,i-k-1$, while it vanishes when $j=i-k,\ldots,i-1$.
 Moreover, $E(Y_{i}^{\prime })=q\Sigma _{r=0}^{k-1}p^{k+r}=p^{k}(
1-p^{k}) $, $i=1,2,\ldots,n-k+1$, whereas $(i\geq 2k)$
\begin{eqnarray*}
E(Y_{j}^{\prime }Y_{i}^{\prime }) &=&E\Biggl( (1-Z_{j-1})\Biggl(
\prod_{l=j}^{j+k-1}Z_{l}+\prod_{l=j}^{j+k}Z_{l}+\cdots+\prod%
_{l=j}^{i-2}Z_{l}\Biggr) Y_{i}^{\prime }\Biggr)
\\
&=& qp^{k}\frac{1-p^{i-j-k}}{1-p}E(Y_{i}^{\prime
})=p^{2k}(1-p^{i-j-k})( 1-p^{k}) ,\qquad  i-2k+1\leq j\leq i-k-1
\end{eqnarray*}
and $E(Y_{j}^{\prime }Y_{i}^{\prime })=0$ for $i-k\leq j\leq i-1$. So, for $%
i\geq 2k$, we get
\begin{eqnarray*}
&&\sum_{j=i-2k+1}^{i-1}\bigl(E(Y_{i}^{\prime }Y_{j}^{\prime
})+E(Y_{i}^{\prime })E(Y_{j}^{\prime })\bigr)
\\
&&\quad =p^{2k}(1-p^{k})\biggl( \biggl(k-1-p\frac{1-p^{k-1}}{1-p}\biggr)+(2k-1)(1-p^{k})%
\biggr) \leq p^{2k}(3k-2)
\end{eqnarray*}
and, thus, for $m=(6k-2)qp^{k}<\log 2$,
\begin{eqnarray} \label{eq35}
&& d_{\mathrm{TV}}\Biggl( \mathcal{L}\sum_{i=1}^{n-k+1}Y_{i}^{\prime
},\mathit{CP}(\lambda ,F_{k})\Biggr)\nonumber
 \\[-8pt]\\[-8pt]
&&\quad\leq \mathit{UB}_{n,k}:=\biggl( 1+\frac{2}{3\lambda }\biggr) (6\lambda kqp^{k})^{2}+%
\frac{{1/4}\Vert \Delta ^{2}f_{\mathit{CP}(\lambda ,F_{k})}\Vert _{1}%
}{1-2(1-\mathrm{e}^{-m})}\frac{\lambda }{q}(6k-3)p^{k},\nonumber
\end{eqnarray}
where $\lambda =(n-k+1)qp^{k}$ and, for $x=1,2,\ldots,k-1$,
\begin{eqnarray*}
F_{k}(x)=P( Y_{i}^{\prime }\leq x |Y_{i}^{\prime }\neq 0)
=P\Biggl( 1+\prod_{j=i+k}^{i+k}Z_{j}+\cdots+\prod_{j=i+k}^{i+2k-2}Z_{j}\leq
x\Biggr) =1-p^{x}
\end{eqnarray*}
is the geometric distribution truncated at $k$ ($F_{k}(k)=1$). It can be
verified that for large $\lambda $, $\mathit{CP}(\lambda ,F_{k})\approx N(\lambda
E(W),\lambda E(W^{2}))$ with $W$ $\sim F_{k}$ and, according to Remark $1,$
we expect that
\begin{eqnarray*}
\bigl\Vert \Delta ^{2}f_{\mathit{CP}(\lambda ,F_{k})}\bigr\Vert _{1}\sim \frac{4}{%
\lambda E(W^{2})\sqrt{2\curpi \mathrm{e}}} \qquad \mbox{as } \lambda \rightarrow \infty .
\end{eqnarray*}

In order to illustrate the above asymptotic relation, we present below a
table with the exact value of the norm $\Vert \Delta ^{2}f_{\mathit{CP}(\lambda
,F_{k})}\Vert _{1}$ and its approximation $4/(\lambda E(W^{2})\sqrt{%
2\curpi \mathrm{e}})$ for several values of $\lambda ,p$ (see Table~\ref{tab1}). We assume that $k\rightarrow
\infty $, that is, $F_{k}$ is the ordinary geometric distribution and thus $%
E(W)=1/q$, $V(W)=p/q^{2}$ and $E(W^{2})=(1+p)/q^{2}$.

\begin{table*}[b]
\tabcolsep=0pt
\caption{}\label{tab1}
\begin{tabular*}{\tablewidth}{@{\extracolsep{4in minus 4in}}lllllllll@{}}
\hline
&${ \lambda =1}$ &  & ${\lambda =5}$ &  & ${ \lambda =10}
$ &  & ${\lambda =100}$ &  \\[-8pt]
&
\multicolumn{2}{l}{\hrulefill}&\multicolumn{2}{l}{\hrulefill}&\multicolumn{2}{l}{\hrulefill}&\multicolumn{2}{l}{\hrulefill}
\\
& \multicolumn{1}{l}{norm} & {approx.} & \multicolumn{1}{l}%
{norm} & {approx.} & \multicolumn{1}{l}{norm} &
{approx.} & \multicolumn{1}{l}{norm} & \multicolumn{1}{l}{approx.}
 \\
 \hline
${  p=0.2}$ & \multicolumn{1}{l}{ 0.97120} & { 0.516204} &
\multicolumn{1}{l}{ 0.115414} & { 0.103241} &
\multicolumn{1}{l}{ 0.054341} & { 0.051620} &
\multicolumn{1}{l}{ 0.005189} & \multicolumn{1}{l}{ 0.005162}
\\
${ p=0.5}$ & \multicolumn{1}{l}{ 1.10364} & { 0.161314} &
\multicolumn{1}{l}{ 0.040737} & { 0.032263} &
\multicolumn{1}{l}{ 0.017866} & { 0.016131} &
\multicolumn{1}{l}{ 0.001628} & \multicolumn{1}{l}{ 0.001613}
\\
${ p=0.8}$ & \multicolumn{1}{l}{ 1.32437} & { 0.021509} &
\multicolumn{1}{l}{ 0.019508} & { 0.004302} &
\multicolumn{1}{l}{ 0.002474} & { 0.002151} &
\multicolumn{1}{l}{ 0.000218} & \multicolumn{1}{l}{ 0.000215}
\\
\hline
\end{tabular*}
\end{table*}

As expected, the above approximation is satisfactory for moderate and large
values of $\lambda $. Moreover, we observe that it becomes better when $p$
decreases. Assuming that $n,k\rightarrow \infty $ with $p\in (0,1)$ fixed,
the compound Poisson approximation error bound in (\ref{eq35}) is of order
\begin{eqnarray*}
\mathit{UB}_{n,k}\sim
\cases{
\dfrac{1}{4}\bigl\Vert \Delta ^{2}f_{\mathit{CP}(\lambda ,F_{k})}\bigr\Vert _{1}%
\dfrac{\lambda }{q}6kp^{k}, &\quad $\mbox{when }\lambda =(n-k+1)qp^{k}\mbox{ is
fixed,}$ \cr
\dfrac{6q}{(1+p)\sqrt{2\curpi \mathrm{e}}}kp^{k}, &\quad $\mbox{when }\lambda \rightarrow
\infty ,\mbox{ such that }\lambda ^{2}kp^{k}\rightarrow 0.$
}
\end{eqnarray*}
For almost the same distance as in (\ref{eq35}), the Stein--Chen method offers a
bound $\mathit{UB}_{CS}$ such that
\begin{eqnarray*}
\mathit{UB}_{CS}\sim \frac{\log ^{+}(\lambda q(1-2p))}{q^{2}(1-2p)}6kp^{k} \qquad
 \mbox{when }p\leq \frac{1}{3}\quad \mbox{or}\quad \mathit{UB}_{CS}\sim \frac{6q}{1-5p}kp^{k}\qquad
   \mbox{when }p\leq \frac{1}{5}
\end{eqnarray*}
(see, for example, Barbour and Chryssaphinou (\citeyear{2001Barbour})). Note that for values of $p>1/3$%
, the Stein--Chen method yields bounds of order $\mathrm{O}(kp^{k}+\mathrm{e}^{-a_{k}\lambda })$
or $\mathrm{O}(\lambda kp^{k})$. The $\mathit{UB}_{n,k}$ is smaller provided that $\lambda
^{2}kp^{k}\approx 0$ and is of order $\mathrm{O}(kp^{k})$ for all values of $p$.

It is worth mentioning that here, instead of Corollary \ref{cor3}(a), we
could employ Corollary \ref{cor3}(b) to obtain a bound even better than $%
\mathit{UB}_{n,k}$. Specifically, it can be proven that for every $1\leq j<i\leq n$,
the random variables $Y_{j}^{\prime }+\cdots+Y_{i-1}^{\prime }$ and $%
Y_{i}^{\prime }$ are NQD. Hence, we can use Corollary \ref{cor3}(b) and,
following an exact parallel to the above procedure, we derive the improved
bound
\begin{eqnarray*}
\mathit{UB}_{n,k}^{\prime }:=12\biggl( 1+\frac{1}{kq}+\frac{2q^{2}}{\lambda }\biggr)
(\lambda kp^{k})^{2}+\frac{{1/2}\Vert \Delta ^{2}f_{\mathit{CP}(\lambda
,F)}\Vert _{1}}{1-2( 1-\mathrm{e}^{-m}) }\biggl( 1+\frac{1+p}{2kq}%
\biggr) \frac{\lambda }{q}kp^{k},
\end{eqnarray*}
which, asymptotically, is about three times smaller than $\mathit{UB}_{n,k}$.

Finally, we can approximate $\sum_{i=1}^{n-k+1}Y_{i}=\sum_{i=1}^{n-k+1}X_{i}$%
, the total number of overlapping success runs within trials $1,2,\ldots,n$, by
$\mathit{CP}(\lambda ,G)$, where $G$ denotes the ordinary geometric distribution with
parameter $p$. In this case, $\mathit{CP}(\lambda ,G)$ is also known as the P\'{o}lya--Aeppli
distribution with parameters $\lambda ,p$ and will be denoted by $\mathit{PA}(\lambda
,p)$. Using the triangle inequality, the distance $d_{\mathrm{TV}}( \mathcal{L}%
\sum_{i=1}^{n-k+1}X_{i},\mathit{PA}(\lambda ,p)) $ is bounded above by
\begin{eqnarray*}
d_{\mathrm{TV}}\Biggl( \mathcal{L}\sum_{i=1}^{n-k+1}Y_{i},\mathcal{L}%
\sum_{i=1}^{n-k+1}Y_{i}^{\prime }\Biggr) +d_{\mathrm{TV}}\Biggl( \mathcal{L}%
\sum_{i=1}^{n-k+1}Y_{i}^{\prime },\mathit{CP}(\lambda ,F_{k})\Biggr) +d_{\mathrm{TV}}(
\mathit{CP}(\lambda ,F_{k}),\mathit{PA}(\lambda ,p)).
\end{eqnarray*}
The first $d_{\mathrm{TV}}$ is bounded by (\ref{eq34}), the second bounded  by (%
\ref{eq35}), whereas for the third, we have ($W_{i},U_{i}$ are independent
random variables with $W_{i}\sim F_{k}$ and $U_{i}\sim G$)
\begin{eqnarray*}
d_{\mathrm{TV}}( \mathit{CP}(\lambda ,F_{k}),\mathit{PA}(\lambda ,p)) =d_{\mathrm{TV}}\Biggl(
\sum_{i=1}^{N}W_{i},\sum_{i=1}^{N}U_{i}\Biggr) \leq \lambda d_{\mathrm{TV}}(
W_{1},U_{1}) =\lambda p^{k}.
\end{eqnarray*}

\section*{Acknowledgments}
The second author was partially supported
by the University of Athens Research
Grant 70/4/8810.

\printhistory

\end{document}